**Chapter 9**
**Christine Proust**[1]

# Early-Dynastic Tables from Southern Mesopotamia, or the Multiple Facets of the Quantification of Surfaces

(preprint)

**Abstract** How were surfaces evaluated before the invention of the sexagesimal place value notation in Mesopotamia? This chapter examines a group of five tablets containing tables for surfaces of squares and rectangles dated to the Early Dynastic period (ca. 2600-2350 BCE) and unearthed in southern Mesopotamia. In order to capture the methods used by ancient scribes to quantify surfaces, special attention is paid to the layout and organization of the tables, as well as to the way in which measurement values are written down. It is argued that these methods vary according to the dimensions of the squares or rectangles concerned: the quantification of small surfaces does not use the same mathematical tools as the quantification of large parcels of land. The chapter shows a reciprocal influence between the metrological systems adopted by the ancient scribes and the methods of calculation of surface they implemented. Some methods may reflect ancient land-surveying practices, and others may testify the emergence of new mathematical concepts applied to all kinds of surfaces, large or small. Ultimately, several different conceptualizations of the notion of surface emerge from the examination of these tables.

## Introduction[2]

This contribution is an attempt to understand how surfaces were quantified in southern Mesopotamia during the Early Dynastic period (ca. 2600-2350 BCE – see chronology in Appendix 9.B), that is, before the invention of the sexagesimal place value notation. The documentation examined is a set of five Early Dynastic tablets considered by historians to be the earliest known mathematical texts.[3] These documents contain tables which establish correspondences between the sides and the surfaces of squares or rectangles. I closely examine the layout and organisation of the tables, as well as the way in which measurement values are written down, in order to capture the methods used by ancient scribes to quantify surfaces. It is argued that these methods vary according to the dimensions of the squares or rectangles: the quantification of small surfaces does not use the same mathematical tools as the quantification of large lands. Several different approaches to measuring surfaces adopted by ancient actors emerge from the examination of these tables. Some may reflect ancient land-surveying practices, and others may testify the emergence of new mathematical concepts applied to all kind of surfaces, large or small.

---

1 Christine Proust, CNRS and Université Paris Diderot - Laboratoire SPHERE UMR 7219 - case 7093 – 5 rue Thomas Mann - 75205 Paris cedex 13, France, email: christine.proust@orange.fr. The research leading to these results has received funding from the European Research Council under the European Union's Seventh Framework Programme (FP7/2007-2013) / ERC Grant Agreement No. 269804.

2 I thank Walter Farber, curator of the cuneiform tablet collection at the Oriental Institute of the University of Chicago, and Joachim Marzahn, curator of the cuneiform tablet collection at the Vorderasiatisches Museum (Berlin), for their kind authorisation to reproduce the photographs. My gratitude goes to the SAW editorial board for their careful reading of the first version of this text, their suggestions and corrections. I specially thank Camille Lecompte for his help in the reading of colophons.

3 Of course, deciding what is the 'earliest mathematical text' is a matter of appreciation on what is a mathematical text. For example, according to the website *Cuneiform Digital Libray Initiative* (http://cdli.ox.ac.uk/wiki/), 'The earliest known mathematical exercise' is a computation of surfaces found in Uruk, dated from the Uruk IV period (ca. 3350-3200 BCE), and now kept at the University of Heidelberg, under the museum number W 19408,76 (+ fragments).

As the arguments are based on the layout and the notations of measurement values adopted in the sources examined, the first section of this chapter specifies the terminology used to describe these elements. The second section summarizes the information available on the five Early Dynastic tables of surfaces known to date, numbered from 1 to 5 (see Table 9.1 below – primary publication and CDLI numbers are provided in the 'list of primary sources' at the end of the chapter). After these preliminaries, the five tables are presented in some detail. These text studies are based on previous work on this documentation, and especially that carried out by Friberg (2007).[4]

| | |
|---|---|
| 1 | VAT 12593 |
| 2 | MS 3047 |
| 3 | Feliu 2012 |
| 4 | A 681 |
| 5 | CUNES 50-08-001 |

**Table 9.1** Sources

First, I analyse the tables organised in a tabular format (Texts 1-3), then the features of the two Early Dynastic IIIb tables, each organised as lists of clauses (Texts 4-5). In conclusion, by contrasting these two groups, and the different tables within each group, I assess the extent to which these sources reflect different methods of computation related to different conceptualisations of surfaces.

## 9.1 On Visual and Textual Aspects of Early Dynastic Tables

### *9.1.1 On Layout*

The texts analysed in this chapter are tables, a very specific genre of text which provides a correspondence between two (or more) sets of values. Here, the sets of values represent the sides of fields and their surface. This correspondence is expressed through a tabular format in Texts 1-3, and through lists of clauses in Texts 4-5.

The texts exhibit significant visual elements such as columns, alignments, horizontal and vertical lines, and headings. The semantics of these elements are part of the meaning of tables in general, but play a particularly important role in archaic documents, where the texts are reduced to very terse clauses.[5] These elements are analysed in the individual text studies below. However, among these elements, the columns deserve special examination in these preliminaries. Two types of columns appear in the texts considered here. The first is the partition of the surface of the written document into vertical strips (see photograph of Text 4, A 681, Fig. 9.13, Sect. 9.4.1). This arrangement, which can be termed 'typographical columns' is common on clay tablets, as well as in other sorts of documents from all periods, such as dictionaries, newspaper and webpages. On clay tablets, the typographical columns run from left to right on the obverse and from right to left on the reverse. The second type of column represents a correspondence between sets of values. This layout, which can be termed 'tabular columns', is specific to tables (see photograph of Text 1, VAT 12593, Fig. 9.6 below). Tabular columns always run from left to right in a given table, regardless of whether the table appears on the obverse or on the reverse. To make a clear distinction between typographical and tabular columns in the following, I number the typographical columns with lowercase figures (i, ii, iii, iv, etc.), and the tabular columns with uppercase figures (I, II, III, IV, etc.).

[4] After I finished writing this chapter, I discovered a very interesting article by Peter Damerow (2016) that raises similar questions to those discussed here, and is based in part on the same texts. In particular, Damerow discusses the problem of the role of multiplication in surface evaluations in the earliest periods of the history of writing (Damerow 2016: 100-106). According to him, 'the evidence on the extent tablets precludes the interpretation' that the areas of the rectangular fields were 'calculated by multiplying length and width' (*ibid.*: 101). Although Damerow relies on observations that are quite similar to those developed in this chapter, some of his premises and conclusions are different. For example, he assumes that 'area measures … depended on length measures from the very beginning' (*ibid.*: 99), while the present chapter, contrasting the treatment of large and small surfaces, tries to historicize the relationship between units of length and surface.

Both types of columns may appear on the same tablet, as shown by the text published by Feliu (2012; see sketch of Text 3, Fig. 9.9 below), where tabular columns are inserted into typographical columns.

### *9.1.2 On Notations of Measurement Values*

As mentioned in the introduction, the notations of measurement values provide key clues on quantification and computation, so this section offers some clarification on these notations. In archaic texts, as in all metrologies, measurement values may be simple or compound.[6] A simple measurement value is written down with only a single piece (e.g. 16 *sar*), and a compound measurement value is written down with several pieces (e.g. 6 *sar* 15 *gin*).[7] Insofar as the order can be detected, the pieces of a compound measurement value are noted, left to right, from the largest to the smallest.

**Components of a simple measurement value (one piece)**

A simple measurement value is generally made up of two components, followed by the name of the quantified quantity or good. For example, on tablet VAT 12593 (Text 1), the first section contains the following signs (Fig. 9.1):

---

5 On the semantics of tables, see Tournès (forthcoming) and notably in this volume Chemla (Chap. 1) and Ossendrijver et al. (Chap. 2); see also Robson (2003).

6 I introduce the distinction between simple and compound measurement values in response to comments by Karine Chemla, who pointed out that the number of pieces used to express a measurement value may provide important clues on computation practices (SAW workshop, January-March 2013). Indeed, this distinction will be helpful in the following.

7 The same measurement value can be expressed in either a simple, or a compound format (for example, in the modern metric system, 2,15 m versus 2 m 1 dm 5 cm). In the same way, the length expressed with the simple measurement value 1×60+20 *ninda* in Ur III documentation (see for example tablet SAT 2 0210 analysed by Stéphanie Rost, Chap. 5), turned out to be expressed as the compound measurement value 1 *uš* 20 *ninda* in Old Babylonian metrological lists. A compound measurement value may also be adopted for the sake of precision (e.g. 2 *sar* 3 *gin* instead of the rounded value 2 *sar* – see for example the texts analysed by Lecompte, Chap. 8 in this volume).

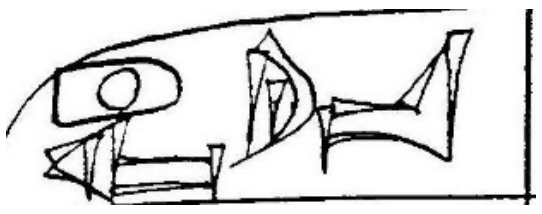

1(geš'u) ninda-DU sag
(600 *ninda* front)
**Fig. 9.1** VAT 12593 (text 1), obv., heading of col. I, copy by Deimel 1923

The two components of this measurement value are the following:

- 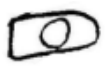 the sign 1(geš'u) represents a numerical value (600),
- 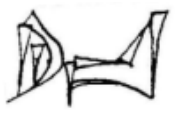 the sign ninda-DU represents a unit of length (about 6 m).

The sign noted under the measurement value specifies that this length refers to the side of a square field:

- 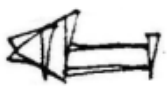 the sign sag, 'front', represents the quantified quantity.

From a graphical point of view, the number is represented by an arithmogram, in curve signs, the measurement unit by a metrogram, in cuneiform writing, and the quantified quantity by a Sumerian logogram, also in cuneiform writing.[8]

This general arrangement can be described as follows:

Numerical value + measurement unit + the quantified quantity or good

Measurement value

The numerical value here is an integer, but in other instances, this value may be noted as a fraction, or as an integer followed by a fraction. In our texts, except for the surfaces of large fields (which are associated with the sign $GAN_2$ – see below), the numbers belong to a non-positional sexagesimal system or 'System S' (Fig. 9.2). 'Non-positional' means that the value of the signs comes only from their shape. The principle of notation is additive: each sign is noted as many times as necessary (e.g. 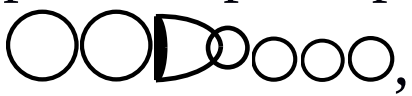, transliterated as 2(šar$_2$) 1(geš'u) 3(u), means 2×3600+1×600+3×10). The system is based on an alternation of factors ten and six, which confers a sexagesimal general structure (hence the name System S given to this numeration by historians, the 'S' meaning 'sexagesimal' somewhat improperly).[9]

| Sign<br>Name | 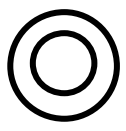<br>1(šar'u) | ×10 ← | 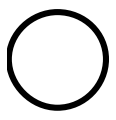<br>1(šar$_2$) | ×6 ← | 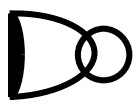<br>1(geš'u) | ×10 ← | 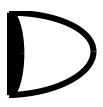<br>1(geš$_2$) | ×6 ← | 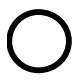<br>1(u) | ×10 ← | 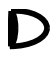<br>1(aš) |
|---|---|---|---|---|---|---|---|---|---|---|---|
| Value | 36 000 | | 3 600 | | 600 | | 60 | | 10 | | 1 |
| In texts | 5 | | 5 | | 3, 5 | | 1, 2, 3, 5 | | 1, 2, 3, 4, 5 | | 1, 2, 3, 4, 5 |

**Fig. 9.2** Systems S

[8] The terms arithmogram and metrogram are imported from Mycenaean studies. The use of this terminology in the analysis of cuneiform writing is attempted, for example, in Proust (2009: Sect. 3.4.5) and Colonna d'Istria (2015).

[9] A discussion on the nature of System S is developed in Ouyang and Proust (forthcoming).

Fractions only appear in Texts 4 and 5 (Fig. 9.3).

| | | |
|---|---|---|
| 1/3 | 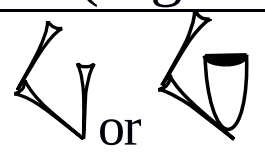 | texts 4, 5 |
| 1/2 | 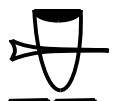 | texts 4, 5 |
| 2/3 | 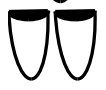 | texts 4, 5 |
| 1/4 (igi-4) | 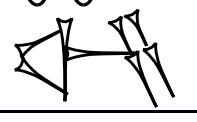 | text 4 |

**Fig. 9.3** Signs for fractions

These signs for fractions turned out to be the ancestors of subsequent cuneiform counterparts, with the notable exception of the sign ‘2(diš)’ for 2/3, whose complex history has been highlighted by Colonna d'Istria (forthcoming).

This general arrangement, with a clear distinction between the numerical and metrological components, is an innovation in the history of writing. In texts predating the Early Dynastic tables, the notation of measurements of length and surface use an integrated system, in which numerical and metrological information is conveyed by the same grapheme.[10] For example, in texts dated from the Uruk IV and III periods, the measures of length are noted with signs belonging to a system similar to System S, without a metrogram. This notation could be considered a defective form of the arrangement (numerical value + measurement unit), but also an integrated system where the signs are arithmo-metrograms. Subsequently, this system was explicitly associated with the measurement unit ‘ninda-DU’ which appears in Early Dynastic tables, and thereafter with ‘ninda’.[11] In the same way, the surfaces are noted with signs belonging to System G (see fig. 4 below) which integrate both numerical and metrological functions. These signs are often accompanied by the sign ‘$GAN_2$’, the ‘pictographic representation of an irrigated field’ (Nissen, Damerow, and Englund 1993: 55), used to indicate that the signs refer to the measurement of fields or lands.

The dissociation between numerical and metrological components as a historical process was brought out by Damerow and Englund (1993: 138):

> [From the Fara period] scribes were becoming increasingly preoccupied with finding better ways of differentiating between numerical signs, metrological signs, and signs for the measured or counted goods.

While in Uruk IV and III texts almost only integrated graphemes are adopted, in our text 1 dated from the Early Dynastic IIIa (or Fara period), as well as in Texts 2 and 3, the measurement units and the quantified quantity or good appear in the headings of the tables (see Fig. 9.1). In Texts 4 and 5, measurement units as such appear in almost all of the clauses. Thus, the five known Early Dynastic tables bear traces of a shift in writing measurement values. However, the function of the sign $GAN_2$ in the quantification of surfaces, and the nature of the so-called ‘System G’ deserve some clarification. This issue is discussed in detail in Appendix 9.A, where I argue that the separation between numerical and metrological components is quite clear in the Early Dynastic tables.[12] Accordingly, the signs of System G are translated as numbers (see Fig. 9.4), and the sign $GAN_2$ is considered, in my commentaries, as a reference field whose surface measures 100 *sar*, that is, a unit of surface defined by: 1 $GAN_2$ = 100 *sar* (or, equivalently, 10 *ninda*-side square).

[10] These integrated graphemes are termed ‘arithmo-metrograms’ in Proust (2009: Sect. 3.4.5) and Colonna d'Istria (2015).

[11] Here I follow the interpretation by Damerow and Englund, who consider that the archaic signs do not represent purely numerical values, but include both numerical and metrological information (Nissen, Damerow, and Englund 1993: 55-58). In this integrated system, ‘1(diš)’ means ‘1 *ninda*’, ‘1(u)’ means ‘10 *ninda*’, and ‘1(geš₂)’ means ‘60 *ninda*’ (Nissen, Damerow, and Englund 1993: 57).

| Sign | 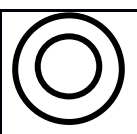 ×10 ← | 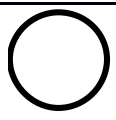 ×6 ← | 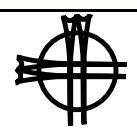 ×10 ← | 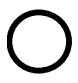 ×3 ← | 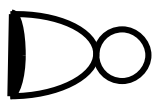 ×6 ← | 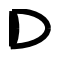 ×2 ← | 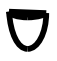 ×2 ← | 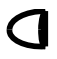 |
|---|---|---|---|---|---|---|---|---|
| Name | 1(šar'u) | 1(šar$_2$) | 1(bur'u) | 1(bur$_3$) | 1(eše$_3$) | 1(iku) | 1(ubu) | ? |
| Value | 10 800 *iku* | 1 080 *iku* | 180 *iku* | 18 *iku* | 6 *iku* | 1 *iku* | ½ *iku* | ¼ *iku* |
| In texts | 5 | 1, 2, 3, 5 | 1, 2, 3, 5 | 1, 2, 3, 5 | 1, 2, 3, 5 | 1, 2, 3, 5 | 3 | 3 |

**Fig. 9.4** System G[13]

### *9.1.3 On Multiplication*

One of the goals of this paper is to capture different methods of quantifying surfaces, some of them involving multiplication. The term 'multiplication' can have several meanings. Jens Høyrup has shown that several notions of multiplication, with different Sumerian / Akkadian names, were used in Old Babylonian mathematical texts (Høyrup 2002: 21-27). In the Early Dynastic texts examined here, we do not find any technical terminology to help us draw such distinctions. However, we shall see that the notation of the measurement values and the organisation of the tables of surfaces clearly reveals different approaches to the operation we uniformly call 'multiplication'.

Assuming Høyrup's distinctions, I use the symbol '×' to denote an arithmetical multiplication, and the symbol '□' (a blank square) to denote a square or a rectangle generated by two linear dimensions. For example, the *sar* measurement unit can be defined as a square whose sides measure 1 *ninda* (ca. 6 m). This relation between length and units of surface can be expressed as

1 *ninda* □ 1 *ninda* = 1 *sar*,

and can be represented as shown in Fig. 9.5.

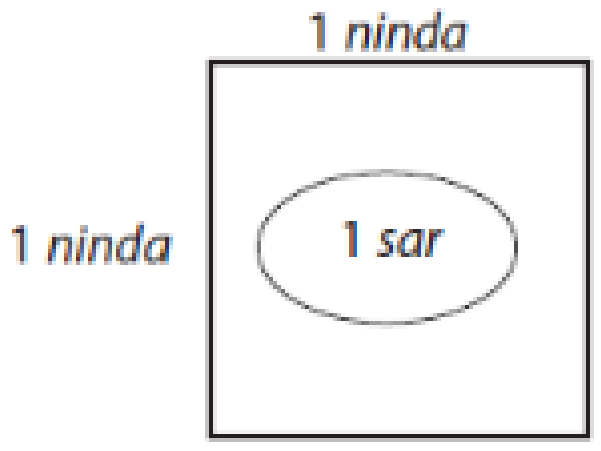

**Fig. 9.5** *ninda* and *sar*

[12] Some historical evolutions should be considered, as the separation between numerical and metrological components is clearer in Text 5 than in Text 1. However, adopting different conventions for representing measurement values written in the same way in different texts would make the reading of this chapter quite uncomfortable.

[13] Note that in Table A of text 5, obv. *ii-v*, multiples of 1(šar$_2$) and 1(šar'u) appear: '*N* gal' means '60 times *N*' and *N* KID' means '60×60 times *N*' (Friberg 2007: 420).

## 9.2 Early Dynastic Tables

Two of the five tablets under study probably date from the Early Dynastic IIIa period (ca. 2600-2500) and come from Šuruppak, while the other three probably date from the Early Dynastic IIIb period (ca. 2500-2340) and come from Adab and perhaps from Zabalam (see Map in Appendix B at the end of the volume). Although their provenance is not always clear, some of them belonging to private collections of unknown origin, the five tablets seem to come from a relatively small geographical area inside the triangle Šuruppak, Nippur and Girsu, in southern Mesopotamia. All of the tables contain a list of measurements of the sides of squares and rectangles, and the measurements of corresponding surfaces. Thus, at first glance, this small *corpus* seems to be quite homogeneous and to reflect a unique culture of computation and quantification. Yet, upon closer inspection, we shall see that these tables bear witness to a dramatic shift in the concept of surface. Clear differentiations are readily apparent when one looks at the layouts and organisations of the tables, as summarised in Table 9.2.

| | **Prove-nance** | **Date** | **Shape of the fields** | **Size order of the entries** | **Format** | **Type of columns** |
|---|---|---|---|---|---|---|
| VAT 12593 (text 1) | Šuruppak | ED IIIa | large square fields | decreasing | Tabular | Tabular columns |
| MS 3047 (text 2) | Šuruppak? | ED IIIa? | large rectangular fields | increasing | Tabular | Tabular columns |
| Feliu 2012 (text 3) | Zabalam? | ED IIIb? | large square and rectangular fields | decreasing | Tabular | Tabular col. inside typographical col. |
| A 681 (text 4) | Adab | ED IIIb | small squares | increasing | List of clauses | Typographical columns |
| CUNES 50-08-001 (text 5) | Zabalam? | ED IIIb? | large and small squares | increasing<br>Set of sub-tables in decreasing size order | List of clauses | Typographical columns |

**Table 9.2** Early Dynastic tables of surfaces (ED = Early Dynastic)

Before entering into a closer examination of the texts, one can observe that they differ from each other in multiple ways. The texts are sometimes organised in a tabular format, sometimes as lists of clauses; the entries are sometimes in decreasing and sometimes in increasing size order; some tables or sub-tables deal with large fields, others with small squares. Moreover, we shall see that the metrological systems are not the same.

## 9.3 Tables of Surfaces in a Tabular Format (Texts 1-3)

The three tables organised in a tabular format provide the measurements of the sides of square and rectangular fields, and the measurements of corresponding surfaces. After a close examination of the layout of the tablets, the organisation of the data, and the notations adopted in the texts, I propose, at the end of this section, hypotheses about the computations of surfaces that these features may reflect.

### *9.3.1 Text 1 (VAT 12593)*

VAT 12593 is an Early Dynastic IIIa tablet from Šuruppak (modern Fara) belonging to the '*Schultexte aus Fara*' corpus published by Deimel (1923: No 82).[14] While no archaeological details were recorded for this tablet by the team in charge of the earliest excavations at Šuruppak, it is highly probable that tablet VAT 12593 was produced by the scholarly milieu of this city who wrote a rich set of lexical lists and Sumerian compositions.[15] This set also includes mathematical texts which show a strong engagement in mathematics by some learned scribes.[16] Šuruppak texts 'reveal fundamental changes in the attitude of the scribe toward arithmetical problems' (Nissen, Damerow and Englund 1993: 137). This scholarly activity can be attributed to a number of identified individuals, whose names appear in the colophon (or subscript) that almost all of the Šuruppak tablets bear.[17] The lexical and mathematical texts from Šuruppak are commonly considered by historians to be 'school texts', as underlined by the title of Deimel's publication. However, as pointed out by Biggs (1974: 29), this characterisation may be misleading insofar as these texts could be considered as a scholarly production rather than young pupils' exercises.

The table in VAT 12593 has been commented abundantly in publications dealing with cuneiform mathematics, and I rely on this previous work in the following presentation.[18] The table is arranged in three tabular columns, which include ten rows on the obverse and six rows on the reverse. Column I provides the lengths of the 'front' (sag) of square fields, column II the lengths of the 'equal side' ($sa_2$) of the square fields, and column III the corresponding surfaces. The sides of the squares are listed in decreasing size order.

---

14 Biggs (1974: 36-39). Biggs listed the 82 Šuruppak tablets published by Deimel containing lexical texts and the present mathematical text (No. 82), along with 24 other literary and lexical texts from Šuruppak published by Jestin (1937 and 1957).

15 The earliest excavations at Fara (ancient Šuruppak) were led by Robert Koldewey and Walter Andrae on behalf of the Deutsche Orient-Gesellschaft in 1902-1903, and subsequently by Erich Smith on behalf of the University Museum, Philadelphia in 1931. German excavators paid attention to the archaeological context of the epigraphic material, but the excavation numbers of some tablets were lost (Krebernik 1998; Robson 2003: 27).

16 'More important, however, than these table texts are a number of others from Fara which tell us something about the nature of mathematical instruction' (Powell 1976: 431).

17 Such a subscript does not appear in tablet VAT 12593 (Text 1), whose reverse is partly broken, but may have existed in antiquity. Other Early Dynastic tables studied here (Texts 2-5) exhibit this kind of subscript.

18 Neugebauer (1935: 91-92); Powell (1976: 430-431); Friberg (1987-1990: 540; 2007: 149); Nissen, Damerow and Englund (1993: 136-139); Robson (2003: 27-30).

## Transliteration of Text 1 (VAT 12593)

Obverse

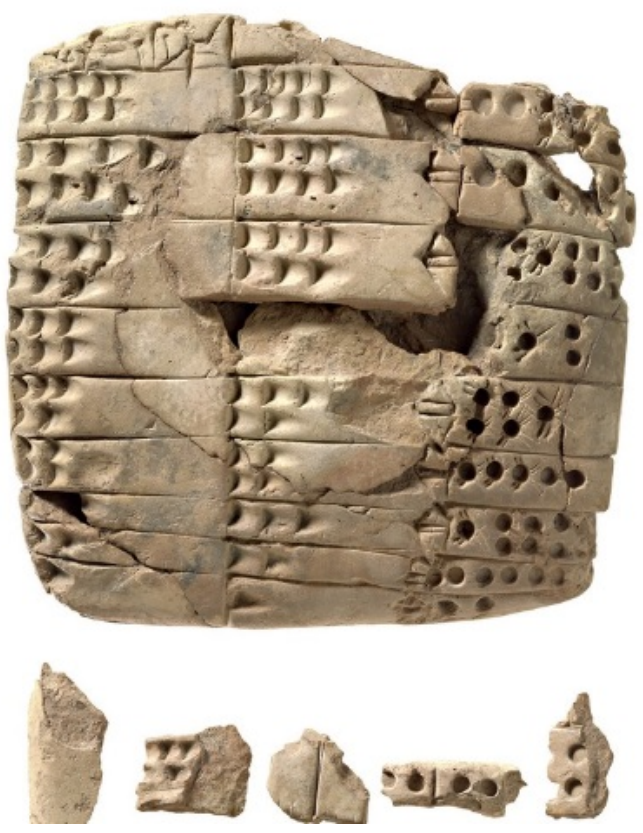

| Col. I | Col. II | Col. III |
|---|---|---|
| 1(geš'u) ninda-DU sag | ⸢1(geš'u)⸣ [sa$_2$] | [3(šar$_2$) 2(bur'u) GAN$_2$] |
| 9(geš$_2$) | 9(geš$_2$) sa$_2$ | 2(šar$_2$) 4 (bur'u) ⸢2(bur$_3$)⸣ |
| 8(geš$_2$) | 8(geš$_2$) sa$_2$ | 2(šar$_2$) ⸢8(bur$_3$)⸣ |
| 7(geš$_2$) | 7(geš$_2$) sa$_2$ | [1(šar$_2$)] 3(bur'u) ⸢8(bur$_3$)⸣ |
| 6(geš$_2$) | [6(geš$_2$) sa$_2$] | 1(šar$_2$) 1(bur'u) 2(bur$_3$) |
| 5(geš$_2$) | 5(geš$_2$) sa$_2$ | 5(bur'u) |
| 4(geš$_2$) | 4(geš$_2$) sa$_2$ | 3(bur'u) 2(bur$_3$) |
| ⸢3(geš$_2$)⸣ | 3(geš$_2$) sa$_2$ | 1(bur'u) 8(bur$_3$) |
| 2(geš$_2$) | 2(geš$_2$) sa$_2$ | ⸢8(bur$_3$)⸣ |
| 1(geš$_2$) | 1(geš$_2$) sa$_2$ | 2(bur$_3$) |

Reverse

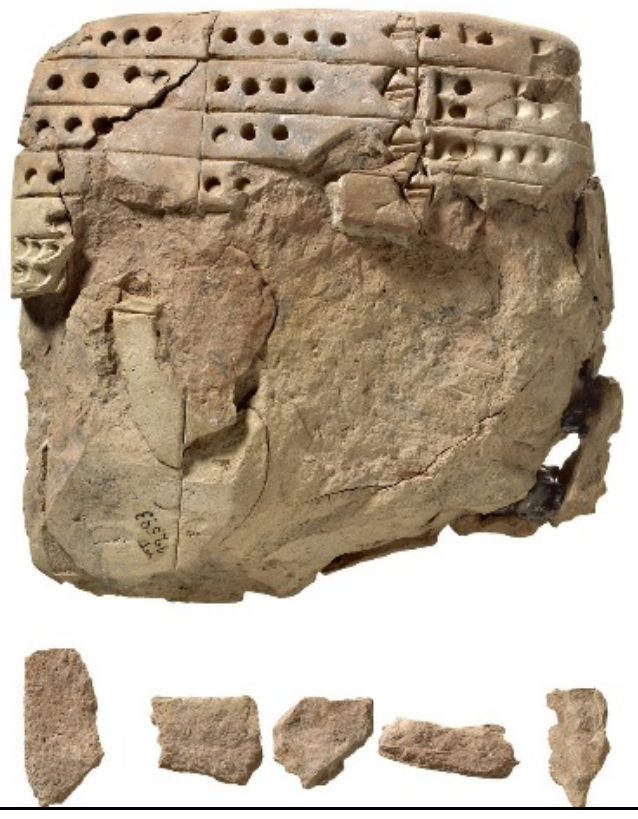

| Col. I | Col. II | Col. III |
|---|---|---|
| 5(u) | 5(u) sa$_2$ | 1(bur$_3$) 1(eše$_3$) 1(iku) |
| 4(u) | 4(u) sa$_2$ | 2(eše$_3$) ⸢4(iku)⸣ |
| 3(u) | 3(u) sa$_2$ | 1(eše$_3$) 3(iku) |
| 2(u) | 2(u) sa$_2$ | [4(iku)] |
| ⸢1(u)⸣ | [1(u) sa$_2$] | [1(iku)] |
| ⸢5(aš)⸣ | [5 sa$_2$] | [1/4(iku)] |

**Fig. 9.6** VAT 12593 (transliteration and photographs C. Proust, courtesy of *Vorderasiatisches* Museum, Berlin)

## Translation of text 1 (VAT 12593)

Obverse

| Col. I | Col. II | Col. III |
|---|---|---|
| 1×600 *ninda* front | 1×600 [equal side] | [3×1080 + 2×180 GAN$_2$] |
| 9×60 | 9×60equal side | 2×1080 + 4×180 + 2×18 |
| 8×60 | 8×60 equal side | 2×1080 + 8×18 |
| 7×60 | 7×60 equal side | [1080 +] 3×180 + 8×18 |
| 6×60 | [6×60 equal side] | 1080 + 180 + 2×18 |
| 5×60 | 5×60 equal side | 5×180 |
| 4×60 | 4×60 equal side | 3×180 + 2×18 |
| ⸢3×60⸣ | 3×60 equal side | 180 + 8×18 |
| 2×60 | 2×60 equal side | ⸢8×18⸣ |
| 1×60 | 1×60 equal side | 2×18 |

Reverse

| Col. I | Col. II | Col. III |
|---|---|---|
| 50 | 50 equal side | 18 + 6 + 1 |
| 40 | 40 equal side | 2× 6 + 4 |
| 30 | 30 equal side | 6 + 3 |
| 20 | 20 equal side | [4] |
| ⸢10⸣ | [10 equal side] | [1] |
| 5 | ⸢5⸣ [equal side] | [1/4] |
| Blank | [Blank] | [Blank] |

The first section of column I contains the full expression of the measurement of length, including the numerical value, the measurement unit and the quantified quantity (600 *ninda* the front, translit. 1(geš'u) ninda-DU sag), while the following sections exhibit only numerical values, namely the number of *ninda*, expressed in System S. The first sections of columns II and III are broken, but may have originally hosted the full expression of the measurement values as the first sections of Text 2 do. The following sections exhibit numerical values in System S with the specification 'equal side' in col. II, and numerical values in System G in col. III. Despite the fact that the fields are square, the sides are designated differently: the first side is termed as 'front' (sag) in the heading, and the second side is termed as 'equal side' (sa$_2$, which means 'equal'). The sign sa$_2$ appears in all of the sections so that a vertical alignment is formed. This asymmetry may reflect the contingent character of the equalities of the sides: Texts 2 and 3 show that this format also applies for rectangular fields.

The metrological system adopted in text 1 is represented in Table 9.3 (the double vertical arrow (↕) indicates the bridge which connects the units of length and surface).

| | |
|---|---|
| Unit of length: | 10 ninda-DU |
| | ↕ |
| Unit of surface: | 1(iku) GAN$_2$ |

**Table 9.3** Measurement units used in Text 1 (VAT 12593)

The sides of the fields vary from 5 *ninda* (30 m) to 600 *ninda* (3600 m): the orders of magnitude of the squares are large lands. Moreover, all the side measures are multiples of 10 *ninda* (except the smallest, 5 *ninda*). This feature shows that the length of 10 *ninda* plays an important role, which also appears in the metrological structure exhibited in Table 9.3. I shall come back to these observations, and their connections with computation, at this end of this section.

### *9.3.2 Text 2 (MS 3047)*

Tablet MS 3047 was published recently by Friberg, who dated the document from the Early Dynastic IIIa period according to the shape of the tablet and the paleography (Friberg 2007: 160). The obverse of the tablet contains a table similar to VAT 12593, but the fields are rectangular instead of square, and the entries are listed in increasing size order instead of decreasing. On the reverse the text is rotated through 90° and is not clearly related to the obverse.

**Transliteration of Text 2 (MS 3047)**[19]

Obverse

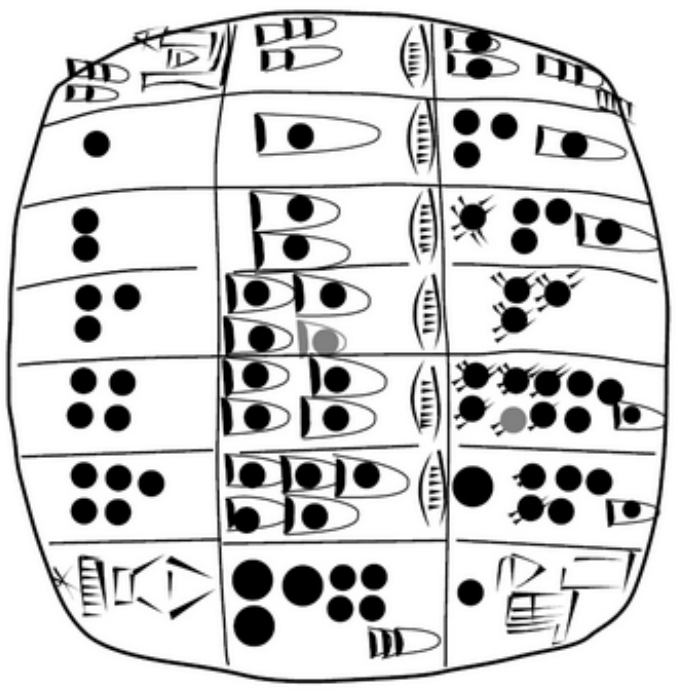

| Col. I | Col. II | Col. III |
|---|---|---|
| sag 5(aš) ninda-DU | 5(ges$_2$) ki | 2(eše$_3$) 3(iku) GAN$_2$ |
| 1(u) | 1(geš'u) ki | 3(bur$_3$) 1(eše$_3$) |
| 2(u) | 2(geš'u) ki | 1(bur'u) 3(bur$_3$) 1(eše$_3$) |
| 3(u) | 3(geš'u) ki | 3(bur'u) |
| 4(u) | 4(geš'u) ki | 5(bur'u) 3(bur$_3$) 1(eše$_3$) |
| 5(u) | 5(geš'u) ki | 1(šar$_2$) 2(bur'u) 3(bur$_3$) 1(eše$_3$) |
| an-se$_3$-gu$_2$ | 3(šar$_2$) 4(bur$_3$) 3(iku) | ŠAR$_2$ GAR GEŠ$^?$ ŠID |

Reverse

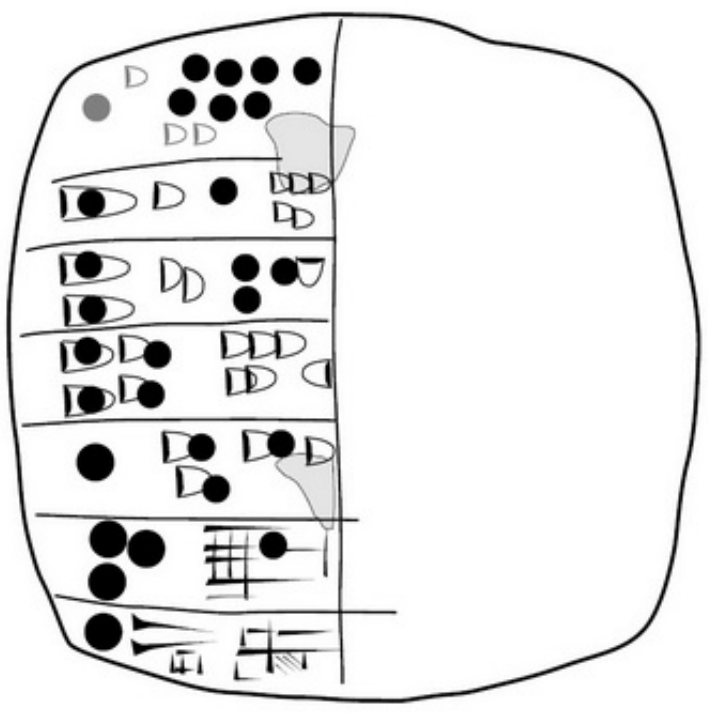

| | |
|---|---|
| […] 7(u) […] | Blank |
| 1(geš'u) 1(geš$_2$) 1(u) 5(aš) | |
| 2(geš'u) 2(geš$_2$) 3(u) 1(ubu)$^?$ | |
| 4(geš'u) 5(geš$_2$) 1/4$^?$ | |
| 1(šar$_2$) 3(geš'u) 1(iku)$^?$ | |
| 3(šar$_2$) E$_2$ 1(u)$^?$ | |
| ŠAR$_2$$^?$ lul$^?$-lul$^?$ | |

**Fig. 9.7** MS 3047 (copy Friberg 2007: 150; reverse rotated 90°)

19 Notes on the transliteration: in obv. III 7, Friberg reads 'x-giš-sanga'. Camille Lecompte observed that the sign 'GIŠ' actually contains some traces and might be identified with another sign, for instance GAN$_2$; the first sign might also not be a numeral but merely šar$_2$; the case could be read ŠAR$_2$? nig$_2$-kas$_7$ GISH+X/GAN$_2$. The reading 'x $^{giš}$nig$_2$-kas$_7$', which may refer to a wooden calculation device ($^{giš}$nig$_2$-kas$_7$) seems to me possible. In rev. 6, Lecompte recognised 'E$_2$' instead of 'KID' suggested by Friberg. In rev. 7, Lecompte cautiously reads lul$^?$-lul$^?$, relying only on a photograph.

**Translation of the obverse of Text 2 (MS 3047)**

| Col. I | Col. II | Col. III |
|---|---|---|
| Front<br>5 *ninda* | 5×60 ground | 2×6 + 3 $GAN_2$ |
| 10 | 1×600 ground | 3×18 + 6 |
| 20 | 2×600 ground | 1×180 + 3×18 + 6 |
| 30 | 3×600 ground | 3×180 |
| 40 | 4×600 ground | 5×180 + 3× 18 + 6 |
| 50 | 5×600 ground | 1×1080 + 2×180 + 3×18 + 6 |
| Total | 3×1080 + 4×18 + 3 ($GAN_2$) | x wood (?) computation |

The metrology and size of the field are the same as in Text 1. The first section of column I contains the full expression of the measurement of length, but in a slightly different order than in Text 1 (front 5 *ninda*, translit. sag 5(aš) ninda-DU), while the following sections exhibit only numerical values, namely the number of *ninda*, expressed in System S. In col. III, the term '$GAN_2$' appears in the heading, thus one may infer that this was also the case in text 1 as supposed in the discussion of this text.

Col. II contains numerical values in System S with the specification 'ground' (ki) where 'equal side' ($sa_2$) is found in Text 1. Friberg (2007: 152) suggests that 'it is possible that ki functions here as a substitute for $aša_5$ [$GAN_2$] "field, surface", and is inserted as a reminder that the product of the length numbers in the first two columns is the surface number in the third column'.) However, it seems to me unlikely that two different signs, ki and $GAN_2$, used in different places – clauses in column II and the heading of column III – have the same function. The function of the term 'ground' (ki) should instead be compared to the function of 'equal side' ($sa_2$) as both are displayed in the same way in Texts 1 and 2: they appear in all of the sections of column II and show a vertical alignment. Thus, $sa_2$ and ki probably have parallel functions in relation to the description of the shape of the fields: square in Text 1, and rectangular in Text 2. The terms 'front' (sag) and 'ground' (ki) may refer to a particular disposition of the fields, as tentatively illustrated by the following diagram (Fig. 9.8).[20]

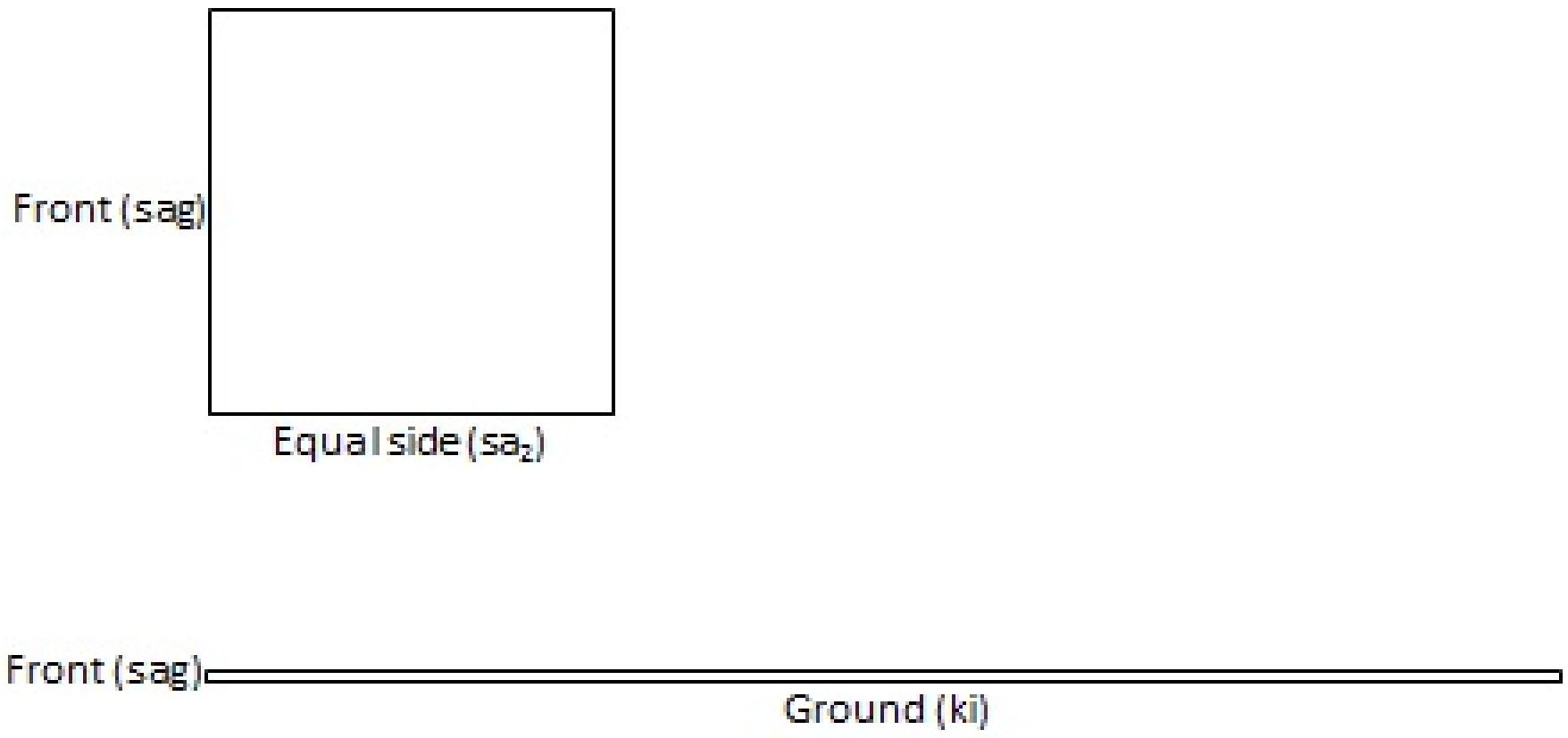

**Fig. 9.8** shape and possible orientation of the fields in text 2

[20] The expression a-ki-ta in the Lagaš texts refers similarly to the orientation of the field (Lecompte, Chap. 8 in this volume).

Note that in Text 3, we again find sa$_2$ and ki, with the same function as in Texts 1 and 2, in relation with square and rectangular fields respectively, which supports Feliu's conclusion: 'we can conclude that sá [sa$_2$] denotes squares, while ki indicates rectangles' (Feliu 2012: 221).

The text on the reverse is not clear. According to Friberg, it seems that System S and System G are mixed in some sections and he observed that most of the values, in System S, form a geometrical progression: the numerical value in line 3 is twice that in line 2, and so on. Even if the exact meaning of the reverse is not clear, the tablet seems to be the result of a mathematical elaboration, underlined by the reference to a 'computation' (nig$_2$-kas$_7$) in the subscript (obv. III 7).

### *9.2.3 Text 3 (Feliu 2012)*

Lluís Feliu has recently published a new table of surfaces belonging to a private collection (Feliu 2012). According to him, a subscript on the obverse mentions the city of Zabalam, where the tablet may have been found during illegal excavations. Feliu dates the tablet from the Early Dynastic IIIb period, like most of the known tablets from Zabalam (Feliu 2012: 219). On the obverse of the tablet there are two typographical columns, each of which contains a table of surfaces in a tabular format. The reverse is mostly destroyed.
Transliteration of the obverse, according to Feliu (2012):[21]

| Column i | | |
|---|---|---|
| Col. I | Col. II | Col. III |
| 1(geš'u) sag GAN$_2$ | 1(geš'u) sa$_2$ | 3(šar$_2$) 2(bur'u) |
| 9(geš$_2$) | 9(geš$_2$) sa$_2$ | 2(šar$_2$) 4 (bur'u) 2(bur$_3$) |
| 8(geš$_2$) | 8(geš$_2$) sa$_2$ | 2(šar$_2$) 8(bur$_3$) |
| 7(geš$_2$) | 7(geš$_2$) sa$_2$ | 1(šar$_2$) 3(bur'u) 8(bur$_3$) |
| 6(geš$_2$) | 6(geš$_2$) sa$_2$ | 1(šar$_2$) 1(bur'u) 2(bur$_3$) |
| 5(geš$_2$) | 5(geš$_2$) sa$_2$ | 5(bur'u) |
| 4(geš$_2$) | 4(geš$_2$) sa$_2$ | 3(bur'u) 2(bur$_3$) |
| 3(geš$_2$) | 3(geš$_2$) sa$_2$ | 1(bur'u) 8(bur$_3$) |
| 2(geš$_2$) | 2(geš$_2$) sa$_2$ | 8(bur$_3$) |
| 1(geš$_2$) | 1(geš$_2$) sa$_2$ | 2(bur$_3$) |
| 5(u) | 5(u) sa$_2$ | 1(bur$_3$) 1(eše$_3$) 1(iku) |
| 4(u) | 4(u) sa$_2$ | 2(eše$_3$) 4(iku) |
| 3(u) | 3(u) sa$_2$ | 1(eše$_3$) 3(iku) |
| 2(u) | 2(u) sa$_2$ | 4(iku) |
| ⸢1(u)⸣ | ⸢1(u) sa$_2$⸣ | 1(iku) |
| [5(aš)] | 5(aš) sa$_2$ | 1/4(iku) |

21 The subscript, quite damaged according to the copy and the photograph, is transliterated by Feliu as '⸢saĝĝa ur⸣-$^{d}$inana saĝĝa zabalam$_4$ (MÙŠ.ZA.UN[UG])' and not translated, except for the name of the city of Zabalam, which, in fact, is not clear.

| Column ii | | |
|---|---|---|
| Col. I | Col. II | Col. III |
| sag<br>1(u) ninda-DU | 1(geš’u) sa$_2$ | [3(bur$_3$) 1(eše$_3$)] |
| 9(aš) | 1(geš’u) sa$_2$ | 3(bur$_3$) |
| 8(aš) | 1(geš’u) sa$_2$ | 2(bur$_3$) 2(eše$_3$) |
| 7(aš) | 1(geš’u) sa$_2$ | 2(bur$_3$) 1(eše$_3$) |
| 6(aš) | 1(geš’u) sa$_2$ | 2(bur$_3$) |
| 5(aš) | 1(geš’u) sa$_2$ | 1(bur$_3$) 2(eše$_3$) |
| 4(aš) | 1(geš’u) sa$_2$ | ⸢1(bur$_3$) 1(eše$_3$)⸣ |
| 3(aš) | 1(geš’u) sa$_2$ | 1(bur$_3$) |
| 2(aš) | 1(geš’u) sa$_2$ | 2(eše$_3$) |
| 1(aš) | 1(geš’u) sa$_2$ | 1(eše$_3$) |
| 3 ur$_2$ hal-la | 1(geš’u) sa$_2$ | 4(iku) 1/2(iku) |
| 2 ur$_2$ hal-la | 1(geš’u) sa$_2$ | 3(iku) |
| 1 ur$_2$ hal-la | 1(geš’u) sa$_2$ | 1(iku) 1/2(iku) |
| 1 kuš$_3$-numun | 1(geš’u) sa$_2$ | 1(iku) |
| ⸢ŠID ur⸣-$^{d}$inana ŠID ⸢zabalam$_4$⸣ | | |

The table in column i, labelled ‘table A’ by Feliu, is a table of surfaces of square fields, in decreasing size order, and is ‘an almost exact duplicate’ of Text 1 (Feliu 2012: 218). Indeed, the terminology, the sizes of the fields, the layouts, the metrological systems, and the notations of the numerical values (i. e. System S in columns I and II, System G in column III) in this table are the same as in Text 1. The fact that the term ‘GAN$_2$’ appears in the heading of column I instead of column III in table A of Text 3 may denote a miscopy. However, Feliu suggests another explanation: ‘As a tentative explanation, it may be that the GÁN sign is used here as a sort of “title” to indicate that the content of the tablet is about “surfaces”’.[22]

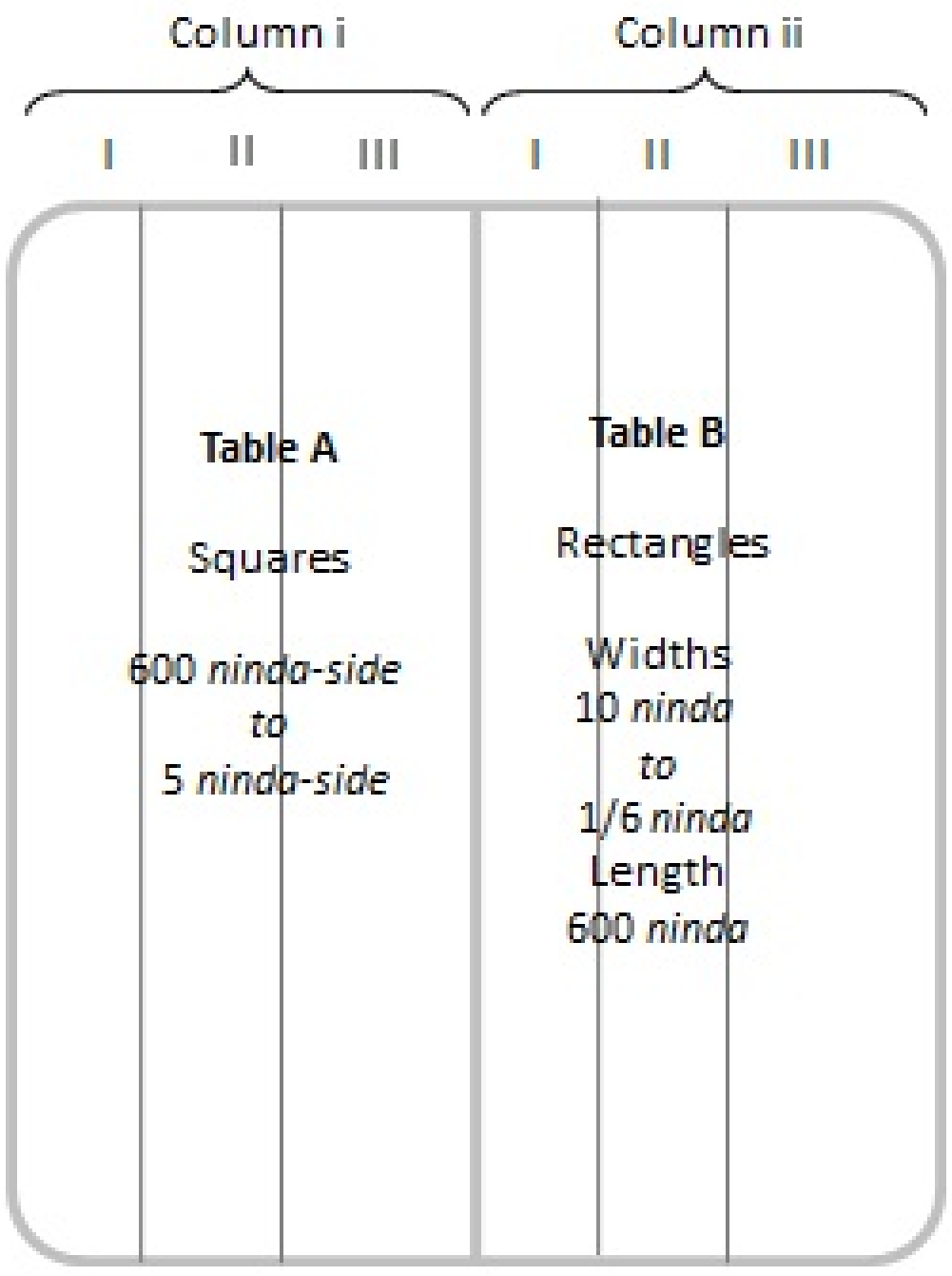

**Fig. 9.9** Typographical columns and tabular columns in Text 3 (obverse)

22 Feliu (2012: 220). Feliu’s explanation fits well with a similar use of GAN$_2$ in Šuruppak tablet Ist Š 188 (TSŠ 188, P010773).

The table in column ii, labelled "table B" by Feliu, is a table of surfaces of rectangular fields similar to Text 2. The same terminology and the same layout are used in Text 2 and table B in Text 3. However, the data differ. The 'front' (sag) measurement of the rectangular fields in table B are listed in decreasing size order and cover a range of smaller values (10 *ninda*, 9 *ninda*, …, 1 *ninda*, and some sub-divisions of *ninda*), and the 'grounds' (ki) are constantly equal to 600 *ninda*. Two new units of length are introduced: ur$_2$ hal-la and kuš$_3$-numun (see Table 9.4).

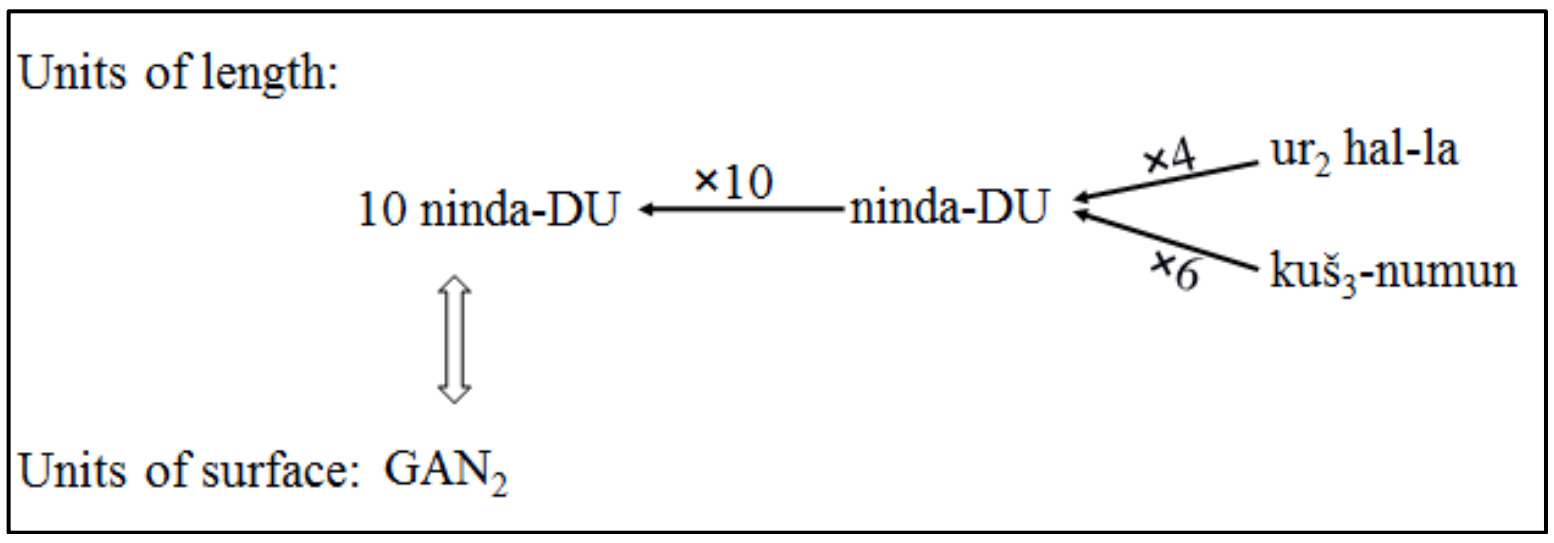

**Table 9.4** Measurement units used in Text 3 (Feliu 2012) and possible bridge

*9.3.4 Computation*

The three texts presented above share many common features, in particular in relation to the way in which measurement values are noted. All of the measurement values are simple (numerical value + measurement unit). The lengths are expressed as numbers of *ninda*, these numbers belonging to System S. The surfaces are expressed as numbers of GAN$_2$, these numbers belonging to System G.

The key to understanding the process for quantifying surfaces is the relationship between the units of length and surface. The definition of the unit of surface '1(iku) (GAN$_2$)' is provided by the three texts (Texts 1, Text 2 obverse, and Text 3 table A), where the item '1(iku)' in col. III corresponds to the entries '1(u)' and '1(u) sa$_2$' in col. I and II. This relationship can be represented as shown in Fig. 9.10, and expressed as:

10 *ninda* □ 10 *ninda* = 1 *iku GAN*

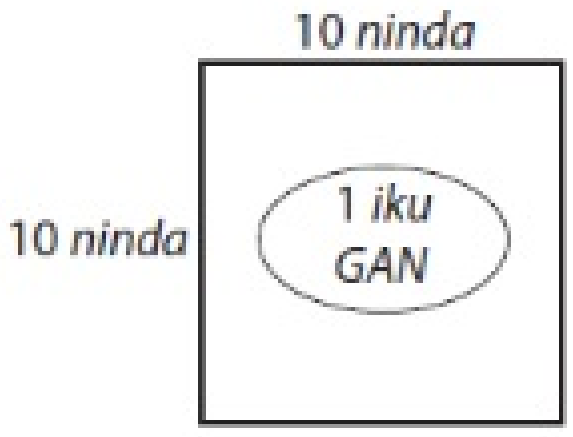

**Fig. 9.10** *ninda* and *GAN*

In later documentation, for example in texts from Early Dynastic Lagaš, a centesimal sub-unit of the *GAN*, the *sar* unit of surface (1 *iku GAN* = 100 *sar*), is used in some texts.[23] This unit of surface is connected to units of length by the fact that 1 *sar* is a 1 *ninda*-side square. However, in Texts 1-3, the *sar* unit of surface does not appear. To what extent was the *sar* used as a surface unit in the context in which table 1 was written, that is, in Šuruppak's scholarly milieus? While the answer to this question is not clear for me, it is certain the *sar* unit of surface appears in Early Dynastic IIIb documentation, for example from Lagaš and Girsu. Thus, this unit was presumably known at least by the authors of Texts 2 and 3.

[23] It is highly probable that the factor 100 between GAN$_2$ and *sar* results from the fact that 1 GAN$_2$ is a 10-*ninda* side square.

From these observations, is it possible to reconstruct the process that was used to quantify the surfaces? A first hypothesis should be that, for each field, the number of *ninda* measuring the 'front' (*N*) was multiplied by the number of *ninda* measuring the 'equal-side' (*M*). The result *N*×*M* (actually, a number of *sar*) would have to be converted into System G. However, there is no trace of the use of the *sar* unit in our texts, nor any trace of such a conversion. Thus, the calculations were probably of a different nature. [24]

Another hypothesis was suggested to me by the reading of an illuminating paper by Duncan Melville about the calculation of rations in Šuruppak texts (Melville 2002). For example, Melville examines a text providing the total quantity of flour to be distributed to 40 men, each man receiving 2 *ban* of flour. He points out that our modern mathematical culture would lead us to imagine that the operation involved is the multiplication of 2 *ban* by 40. However, this multiplication produces a result which does not exist in ancient metrological systems. Melville argues that in fact, solving such a problem does not require a multiplication, and that repeated additions easily lead to the result, requiring only attested measurement values. Following Melville's approach, we can examine the possible methods used for finding the surfaces provided by the Šuruppak table.

As highlighted above, the key issue is the relationship between length and surface measurement values, hence, between systems S and G. A striking feature of the table in Text 1 (as well as of the tables in Texts 2-3) is that there is no item for a square of 1 *ninda*-side and, correspondingly, no trace of the *sar* unit. The bridge between length and surface is not the 1 *ninda*-side square, as usually assumed, but the 10 *ninda*-side square whose surface is 1 *iku GAN*.[25]

If we consider the 10 *ninda*-side square to be pivotal in the relationship between lengths and surface measurements values, it is possible to understand the articulation between Systems S and G. Indeed, if we represent each element of the unit of surfaces (surfaces corresponding respectively to 1 *iku GAN*, 1 *eše GAN*, 1 *bur GAN*, 1 *bur'u GAN* and 1 *šar GAN*), we can observe that each of these elements results from quite simple combinations of the lower elements (see Fig. 9.11, where I introduced a change of scale for the surface 1 *bur GAN*, in order to be able to represent all of the units on the same page). [26] Of course, several different shapes are possible for the elementary surfaces in Fig. 9.11. For example, the sides of the rectangle whose surface is 1 *bur GAN*, measure 60 *ninda* and 30 *ninda* in Fig. 9.11, but the sides of this rectangle could also be 10×60 *ninda* and 3 *ninda*. It is perhaps not by chance that this reshaped rectangle appears in Text 3, as we shall see below.

---

24 The following search for another reconstruction of calculations in Texts 1-3 was provoked by the scepticism of my SAW colleagues faced with a fist version of this paper. As pointed out above, this scepticism was expressed by Damerow (2016: 101).

25 In connection with the essential role of the 10 *ninda*-side square in texts 1-3, it is meaningful to underline that, in administrative texts, the length 10 *ninda* is sometimes noted with the unit *eš* (translit. $\text{eš}_2$). For example in Lagaš metrology, the 1 *eš*-side square provides the bridge between units of length and surface (1 *eš* □ 1 *eš* = 1 $\text{GAN}_2$) (see Lecompte, Chap. 8 in this volume).

26 Similar diagrams appear in Powell (1987-1990: 479-481), who systematically provides geometrical representations of the relationships between length and surface units attested in Mesopotamia, and by Friberg (1997-1998: 6, Fig. 1.2), who provides a very nice synthetic representation of the basic relationships between units of surfaces.

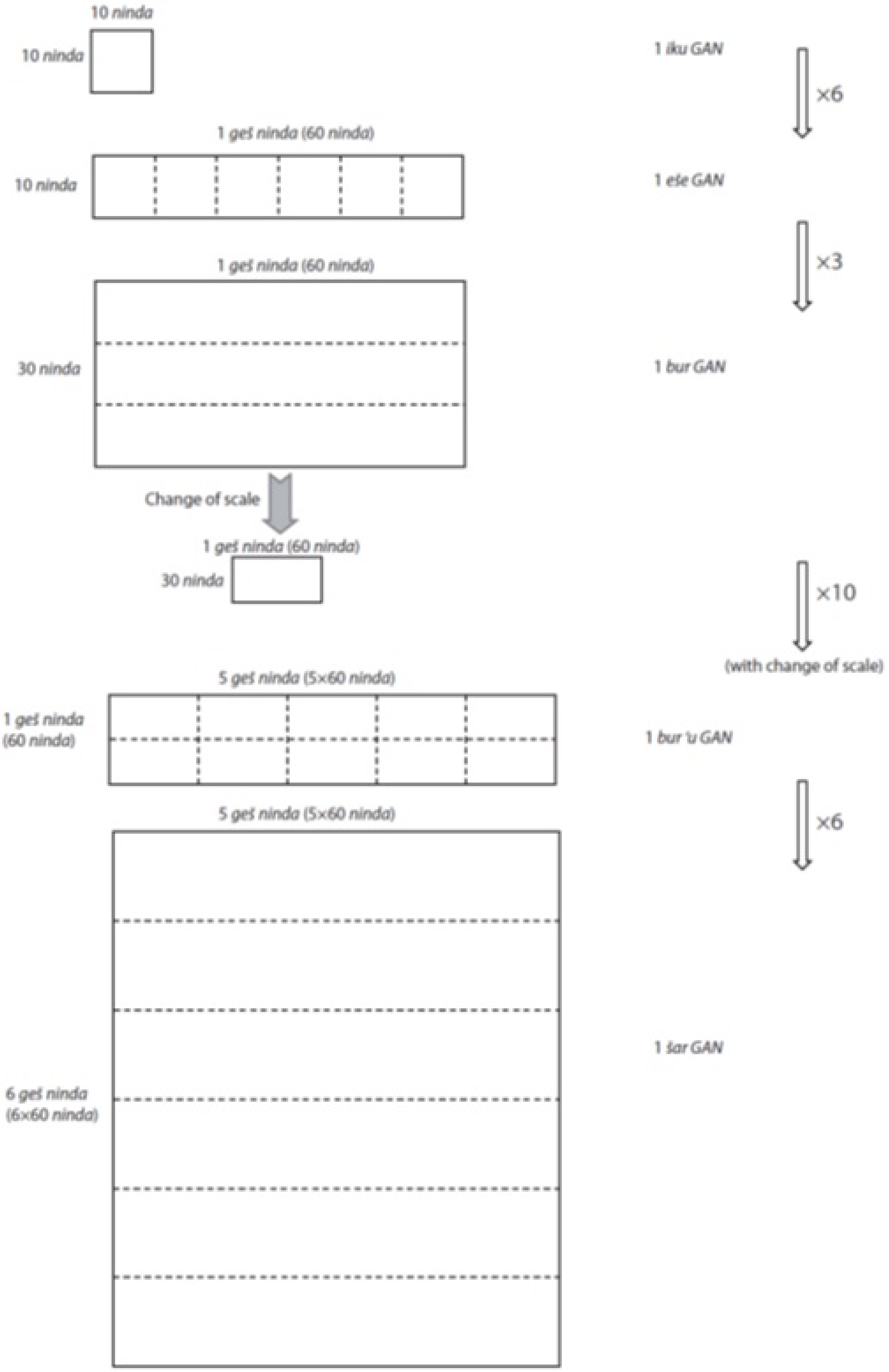

**Fig. 9.11** Relationships between length and surface elements used in Text 1

With these elementary squares and rectangles in mind, it becomes possible to evaluate the surfaces in Text 1 without any calculation. For example, each surface can be obtained by bordering the previous one by the relevant band, as shown in Fig. 9.12a-b (the band is 10

*ninda*-width in Fig. 9.12a, and 60 *ninda*-width in Fig. 9.12b). In this interpretation, the surfaces are constructed from the smaller to the larger by adding a band. In this case, the order of calculation would be reverse of the order of writing adopted in the table.

| | **Side** | **Surface** |
|---|---|---|
| 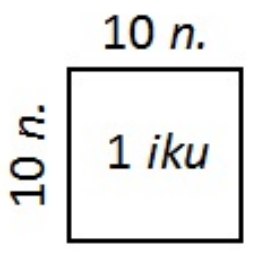  | 10 *ninda* | 1 *iku GAN* |
| 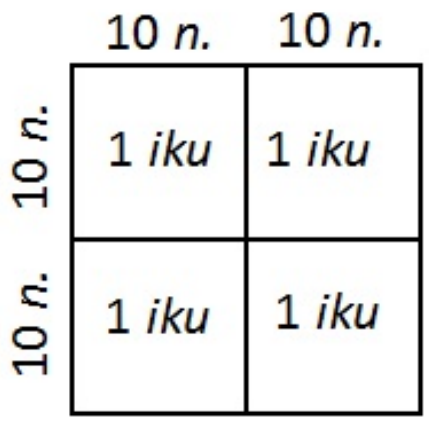  | 20 *ninda* | 4 *iku GAN* |
| 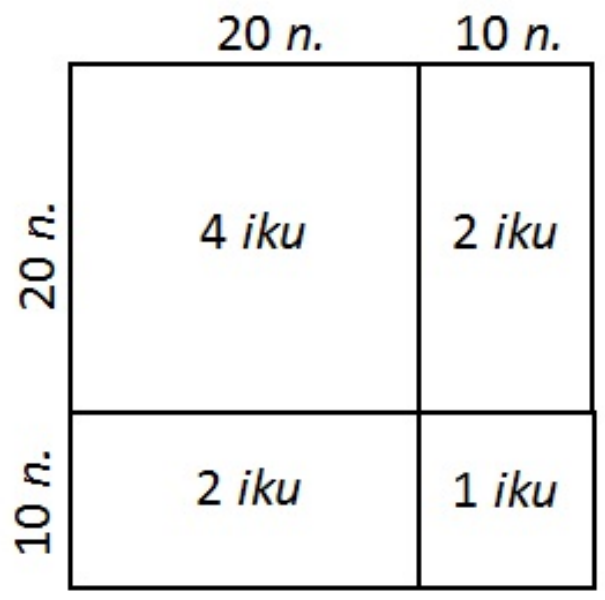  | 30 *ninda* | 1 *eše* 3 *iku GAN* |
| 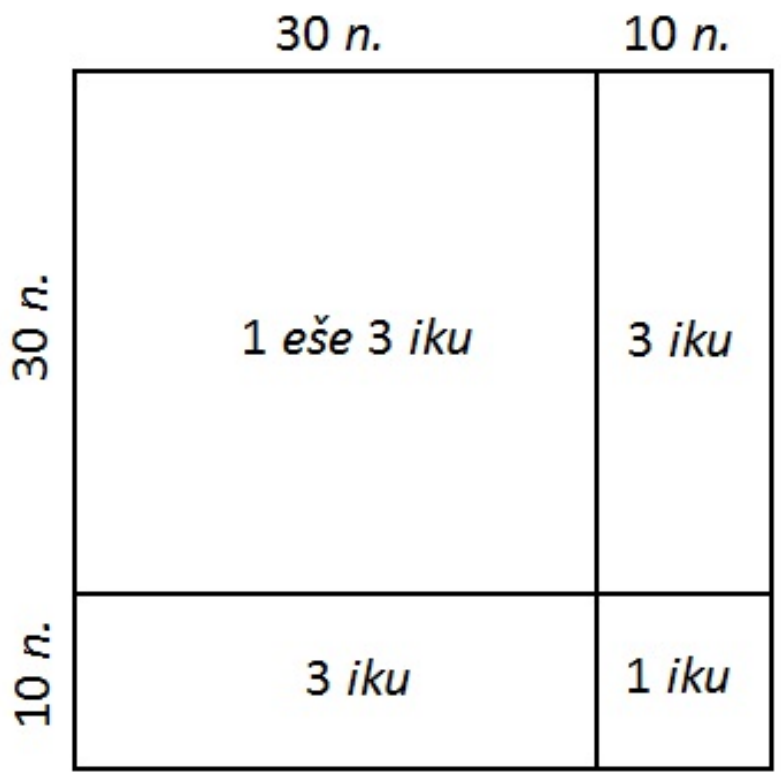  | 40 *ninda* | 2 *eše* 4 *iku GAN* |
| Etc. | | |

**Fig. 9.12a** Representation of square fields tabulated in Text 1 (sides from 10 *ninda* to 40 *ninda*)

<table>
<tr><th></th><th>Side</th><th>Surface</th></tr>
<tr><td>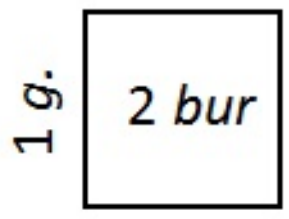
</td><td>60 ninda =<br>1 geš ninda</td><td>2 bur GAN</td></tr>
<tr><td>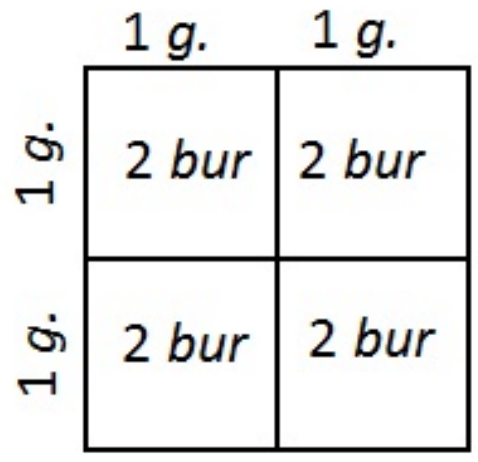
</td><td>2 geš ninda</td><td>8 bur GAN=</td></tr>
<tr><td>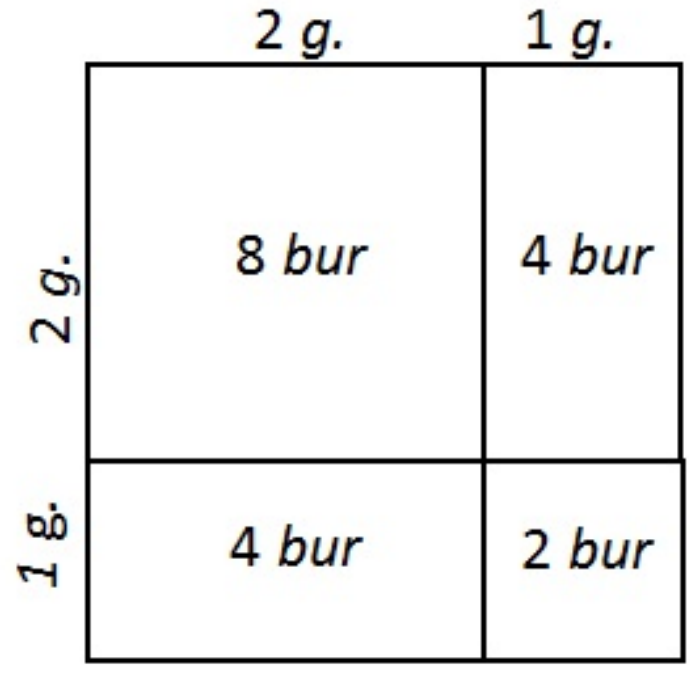
</td><td>3 geš ninda</td><td>1 bur'u 8 bur GAN</td></tr>
<tr><td>Etc.</td><td></td><td></td></tr>
</table>

**Fig. 9.12b** Representation of square fields tabulated in Text 1 (sides from 60 *ninda* to 180 *ninda*)

The difficulty of *calculating* the surface using multiplication comes from the fact that the factors defining each element of System S (respectively 10, 6, 10, 6, 10 – see Fig. 9.2 above) does not seem to be correlated with the factors defining each element of System G (respectively 6, 3, 10, 6, 10 – see Fig. 9.4 above). However, if we adopt a geometrical instead of an arithmetical point of view, these factors produce elements of surfaces that are easier to handle. In addition to the basic correspondence of a 10 *ninda*-side square with 1 *iku GAN*, another simple bridge appears: a 60 *ninda*-side square corresponds to a 2 *bur GAN* surface.

Text 3 contains two tables side by side, table A being identical to the table in Text 1, and table B providing surfaces of rectangles. The striking point is that the rectangles in table B include, among others, surface elements measuring 1 *iku GAN*, 1 *eše GAN*, and 1 *bur GAN*, that is basic elements of System G. Thus, we can hypothesize that table B (or another similar one) may have served as a tool for evaluating surfaces in table A, and, in this case, that the generation of table A in Text 3 was similar to the generation of the table in Text 1. However, other hypotheses are possible. Table A in Text 3 may have been reproduced by copying (or the memorisation of) an earlier table, and not generated by the author of Text 3. Or, as Text 3 possibly comes from Zabalam, and may date from the Early Dynastic IIIb period, it is possible that the methods used by its author were different from those current some decades earlier in Šuruppak.

### *9.3.5 Size of the Fields*

The smallest surface found in the three tables is 1/4 *iku GAN* (ca. 900 m²). These orders of magnitude relate to large lands, and not to plots or gardens cultivated by individuals. This is the scale of the activities of land-surveyors attached to the temple or the central administration, such as the officials of Lagaš described by Lecompte (the sides of these fields cover a range from 3 *ninda* to 100 *ninda*). Most of the rectangular fields have apparently unrealistic shapes: in Text 2, the lengths are exactly sixty times longer than the widths, and in Text 3, table B, the lengths are more than sixty times longer than the widths (in the last item, the length is 3600 times the width). These rectangles are much more elongated than the ‘strip’ fields found in ED IIIb and Ur III documentations, where the length / width ratios do not exceed twenty (see Lecompte, Chap. 8, who refers to Liverani 1996: 15). Thus, these very elongated rectangles are not to be considered as fields, but rather as pieces of surface used for evaluating large lands.

What was the function of texts 1-3? The texts denote a mathematical elaboration, but at the same time, the convenient tabular layout suggests that the tables served as tools for quantifying large lands, saving the land-surveyors tedious computations. Anyway, the three texts probably have different dates and different provenances, and may have been produced or used for slightly different purposes. A table such as Text 1, of which a duplicate appears in Text 3, seems to have been copied and transmitted from one generation to another. This table may have been composed originally as a mathematical elaboration, in the context of the scholarly milieu of Šuruppak, and subsequently used as a tool by land-surveyors, or a similar scenario suggested by Damerow and Englund (1993: 139):

> The exact purpose of this table of areas fields is not known. We may exclude the possibility that it served as some sort of table of calculations used to consult particular values. The list was more likely to have been written as an exercise containing easily determinable field surfaces every land-surveyor was required to know which could be added together in calculating complicated surfaces.

## 9.4 Tables of Surfaces as Lists of Clauses (Texts 4 and 5)

In their content, Texts 4 and 5 are tables of surfaces similar to those examined above. However, they differ completely if we consider the layout, the notations and the range of data. By contrasting the features of these two texts, in comparison with the previously examined tabular texts, I try to capture different conceptualisations of surfaces.

### *9.4.1 Text 4 (A 681)*

Tablet A 681 comes from Adab and is dated to the Early Dynastic IIIb period. Like Text 1, the table is well known to specialists and often cited.[27] However, to my knowledge, details of the notations have not been examined before. The text contains a list of clauses successively providing the sides of squares and the corresponding surface in increasing size order. The list is arranged in typographical columns, three on the obverse, and three on the reverse, of which one is blank.

[27] The publication history, summarised by Feliu (2012: 218), includes Luckenbill (1930 : 70 copy); Edzard (1969: 101-104 edition and study); Friberg (2005: 19 copy and study); Friberg (2007a: 357-360 copy, edition and study); Robson (2007 transliteration, translation and photograph).

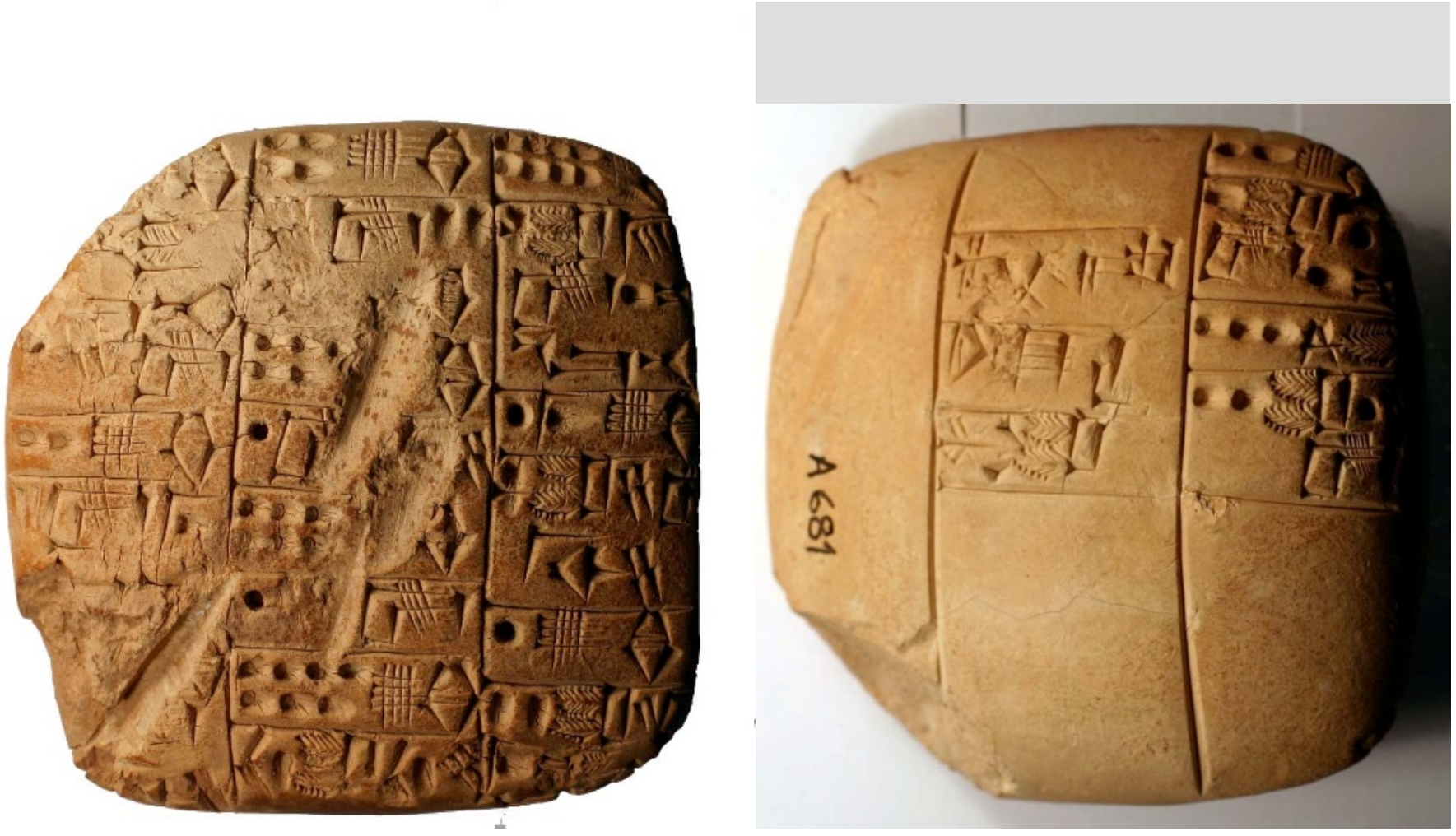

**Fig. 9.13** A 681(photographs by C. Proust, courtesy of Oriental Institute of the University of Chicago)

## Transliteration of Text 4 (A 681)

Obverse

| Column i | Column ii | Column iii |
|---|---|---|
| 1. [1(aš)] kuš$_3$# sa$_2$ | 7. 4(aš) kuš$_3$ sa$_2$ | 15. 8(aš) kuš$_3$ sa$_2$ |
| 2. [1(aš)] sa$_{10}$-ma-na 1(u) 5(diš) gin$_2^{sic}$ (še$^!$) | 8. 7(diš)* gin$_2$ 2(diš) [sa$_{10}$-ma]-na | 16. 1/2 sar la$_2$ 3(diš) gin$_2$ 1(aš) sa$_{10}$-ma-na |
| 3. 2(aš) kuš$_3$ sa$_2$ | 9. 5(aš) [kuš$_3$] sa$_2$ | 17. 1(u) la$_2$ 1(aš) kuš$_3$ sa$_2$ |
| 4. 2(diš) gin$_2$ la$_2$ 1(aš) sa$_{10}$#-ma#-na# | 10. 1(u) gin$_2$# [...]5 | 18. 1/2 sar 4(diš) gin$_2$ la$_2$ igi-4 |
| 5. [3(aš) kuš$_3$] sa$_2$ | 11. 6(aš) [kuš$_3$] sa$_2$ | 19. 1(u) kuš$_3$ sa$_2$ |
| 6. [4(diš) gin$_2$] la$_2$ igi-4(diš) <gin$_2$> | 12. 1(u) [5(diš)] gin$_2$ | 20. 2(diš) (=2/3) sar 2(diš) gin$_2$ la$_2$ 1(aš) sa$_{10}$-ma-na |
| | 13. 7(aš) kuš$_3$ sa$_2$ | |
| | 14. 1/3 sar 1/3 <gin$_2$> 5(diš) <še> | |

Reverse

| | Column ii | Column i |
|---|---|---|
| | 25. nam-mah | 21. 1(u) 1(aš) kuš$_3$ sa$_2$ |
| | 26. nig$_2$-kas$_7$ | 22. 1(aš) sar la$_2$ 1(u) gin$_2$ 1(aš) sa$_{10}$-<ma-na> 1(u) 5(diš) <še> |
| | 27. mu-sar | 23. 3(aš) gi sa$_2$ |
| | | 24. 2(aš) sar 1(u) 5(diš) gin$_2$ |

## Translation of Text 4 (A 681)

Obverse

Column i

1. [1] *kuš* equal side,
2. [1] *samana* 15 *še*$^!$.
3. 2 *kuš* equal side,
4. 2 *gin* minus 1 *samana*.
5. [3 *kuš*] equal side,
6. [4 *gin*] minus 1/4 (*gin*).

Column ii

7. 4 *kuš* equal side,
8. 7$^{sic}$ (6$^!$) *gin* 2 *samana.*
9. 5 [*kuš*] equal side,
10. 10 *gin* [1 *samana* 15 *še*].
11. 6 [*kuš*] equal side,

12. 1[5] *gin*.
13. 7 *kuš* equal side,
14. 1/3 *sar* 1/3 <*gin*> 5 <*še*>.

Column iii

15. 8 *kuš* equal side,
16. 1/2 *sar* minus 3 *gin* 1 *samana*.
17. 10 minus 1 *kuš* equal side,
18. 1/2 *sar* 4 *gin* minus igi-4.
19. 10 *kuš* equal side,
20. 2/3 *sar* 2 *gin* minus 1 *samana*.

Reverse

Column i

21. 11 *kuš* equal side,
22. 1 *sar* minus 10 *gin* 1 *samana* 15 (*še*).
23. 3 *gi* equal side,
24. 2 *sar* 15 *gin*.

Column ii

25. Nam-Mah wrote the computation.

Besides the non-tabular format already underlined, two other features distinguish this text from the previous ones. The formula describing the sides of the squares is slightly different: numerical value + measurement unit + 'equal side' (eg. 4 *kuš* equal side, translit. 4(aš) kuš$_3$ sa$_2$). Here, only one side is given, and this side is labelled 'equal side' (sa$_2$); the term 'sa$_2$' is used here specifically for squares, as in subsequent documentation. More importantly, the sides are smaller than 1 *ninda*, except for the last (3 *gi* = 1 ½ *ninda*), and the surfaces are often expressed as compound measurement values, using subtractive notations. The units of length adopted are sub-units of the *ninda* (2 *gi* = 1 *ninda* and 6 *kuš* = 1 *gi*), and the units of surface are *sar*, completed by three sub-units (60 *gin* = 1 *sar*, 3 *samana* = 1 *gin*, and 60 *še* = 1 *samana*), as shown in Table 9.5.[28]

| | |
|---|---|
| Units of length: | (ninda ←×2)− gi ←×6− kuš$_3$ |
| | ↕ |
| Units of surface: | sar ←×60− gin$_2$ ←×3− sa$_{10}$-ma-na ←×60− še |

**Table 9.5** Measurement units used in Text 4 (A 681)

As the factors defining the units of length are 2 and 6, and the factors defining the units of surface are 60, 3 and 60, the units of surface are not squares of units of length. Thus, the only bridge between the measurements of length and surface is the relation between the largest units, *ninda* and *sar* (1 *ninda* □ 1 *ninda* = 1 *sar*). Actually, the entry 1 *ninda* is missing: the text jumps from the entry 11 *kuš* (rev. col. i 21) to the entry 3 *gi* (rev. col. i 23), omitting the entry1 *ninda* (12 *kuš*), surface 1 *sar*, which should appear quite naturally after entry 11 *kuš*. I shall come back to the meaning of this omission.

---

28 In comparison with Old Babylonian metrology, only two slight differences can be noted: the unit *gi* between *ninda* and *kuš*, and the unit *samana* between *gin* and *še* disappeared in Old Babylonian metrological lists (but the unit *gi* is still used in OB documents). The meaning of sa$_{10}$-ma-na, word for word 'exchange *mana*', is discussed in Friberg (2007: 358).

### *9.4.2 Computation*

Friberg (2007: 358-359) supposed that the lengths of the sides were converted into fractions of *ninda*, and that the computations were carried out on the fractions. However, such computation on fractions did not leave any trace in this text. The puzzling trace we have is the use of subtractive notation for the surfaces.[29] For example, the surface of a 2 *kuš*-side square (obv. i, 3) is expressed as '2 *gin* – 1 *samana*' instead of '1*gin* + 2 *samana*'. These subtractive notations for measurement values do not appear in tables of large fields and may be related to the specific problem of evaluating small surfaces.

The subtractive notation is hard to understand if the surfaces result from a multiplication or computations on fractions. The surfaces seem rather to result from cutting and pasting of elementary squares and rectangles. For example, let us consider the 11 *kuš*-side square (rev. i, 21-22). The surface of the 11 *kuš*-side square can be obtained from the surface of the 12 *kuš*-side square (that is, 1 *ninda*-side square, whose surface is 1 *sar*) by the following cut-and-paste process.[30]

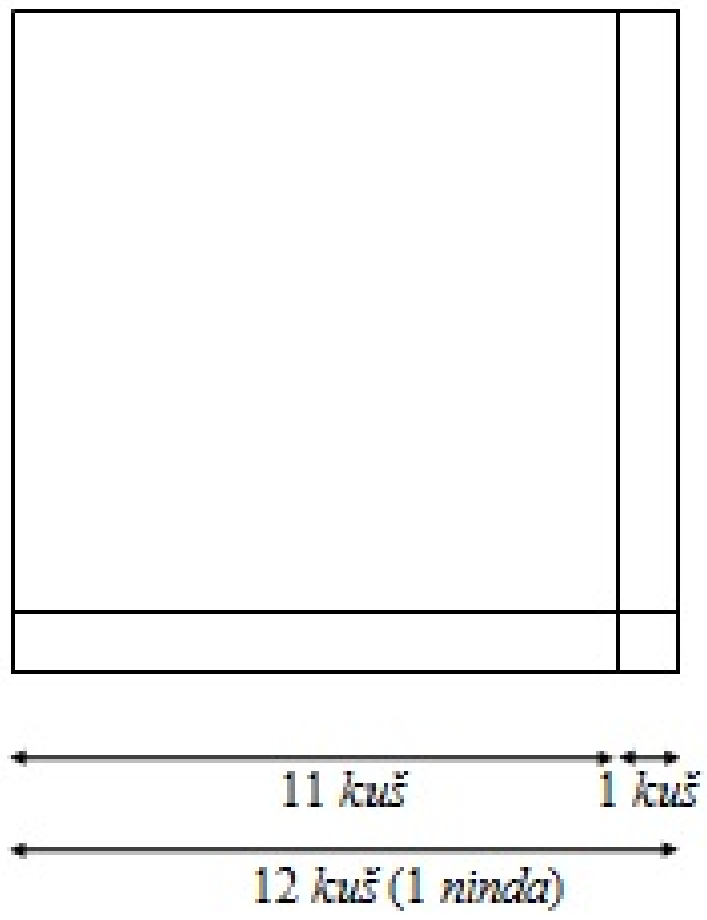

**Fig. 9.14** Surface of the 11 *kuš*-side square

The surface of the 11 *kuš*-side is equal to the surface of the 1 *ninda*-side square *minus* the surface of the edge (see Fig. 9.14). This edge is composed of two strips whose sides are 1 *kuš* and 1 *ninda* (12 *kuš*). The surface of the 1 *kuš*-square is provided by the first clause of the text, which states that the surface of the 1 *kuš*-square is 1 *samana* 15 *še*. This statement is not trivial and seems to have been based on previous knowledge. The surface of each strip (1 *kuš* wide and 1 *ninda* long) is one-twelfth of a *sar*, that is, 5 *gin* (see Fig. 9.15-16 below). Here again this result is not trivial and seems to have been based on previous knowledge.

---

[29] I distinguish the subtractive notation for measurement values (i.e. 1/2 sar minus 3 $gin_2$), from the subtractive notation for numerical values (i.e. 10 minus 1) which is very common in cuneiform texts. Subtractive notations for surfaces appear in A 681 obv. i, 4, 6; iii, 16, 17, 18, 20, and subtractive notations for numbers appear in A 681, obv. i, 6 (9 = 10-1) and 13 (49 = 50-1).

[30] This idea comes from discussions in the SAW workshop in 2013, including suggestions by Carlos Gonçalves.

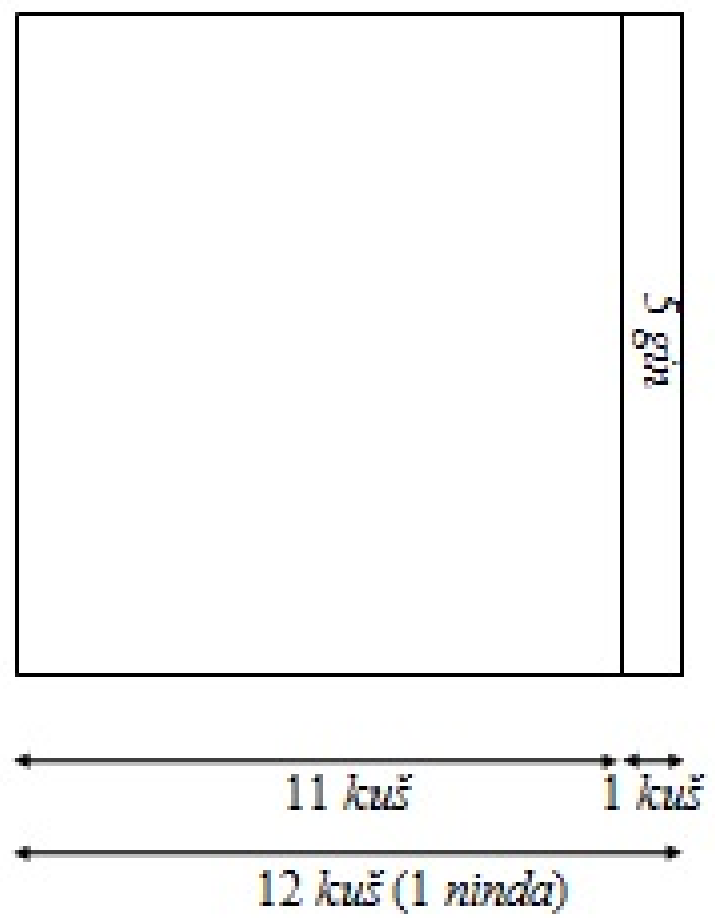

**Fig. 9.15** Strip 1 *kuš* width and 1 *ninda* (12 *kuš*) long

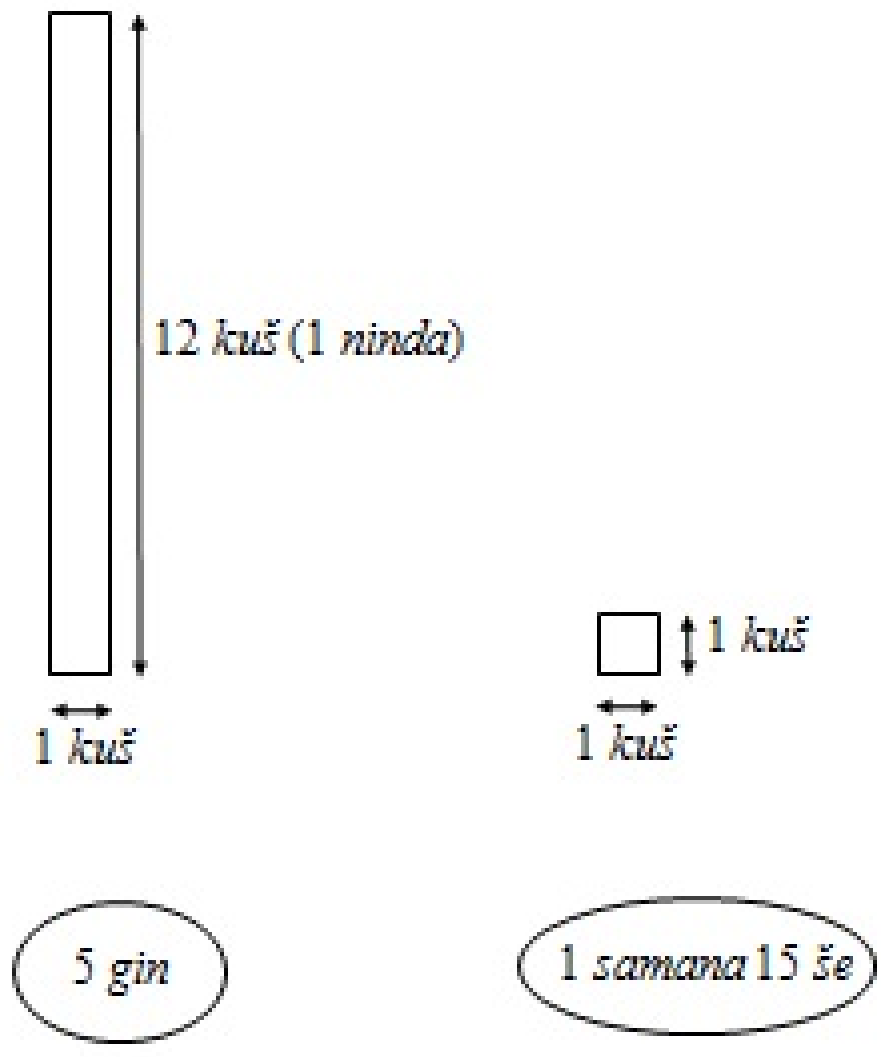

**Fig. 9.16** Elementary surfaces

The two strips of 5 *gin*-surface are cut from the large 1 *sar*-surface square, but this operation causes a 1 *kuš*-side small square (surface 1 *samana* 15 *še*) to be removed twice; thus, this small square must be added back (see Fig. 9.15-16).

The measure of the surface found in the tablet is '1 *sar* – 10 *gin* + (1 *samana* 15 *še*)',[31] that is, 1 *sar* from which 10 *gin* (twice 5 *gin*) is subtracted and 1 *samana* 15 *še* is added. This notation exactly represents the geometrical process described above.

Another example is the 8 *kuš*-side square, whose surface, according to the tablet (obv. iii, 15-16), measures '1/2 *sar* – (3 *gin* 1 *samana*)'. A 2 *kuš* wide strip is cut from the 8 *kuš*-side square and pasted in order to compose a rectangle whose sides measure 6 *kuš* (½ *ninda*) and

---

31 Note that the subtractive notation may be ambiguous. Here, the smallest piece of the measurement value, '1 samana 15 še', is to be added, and not subtracted, to the principal piece, '1 sar'. The subtraction 'la$_2$' acts only on '10 gin'. This contrasts with the usual situation, where subtraction 'la$_2$' acts on all of the pieces that follow it. To avoid ambiguity, I introduce parenthesis to indicate the piece to be subtracted from 1 *sar*, and the pieces to be added.

12 *kuš* (1 *ninda*), and surface 1/2 *sar* (Fig. 9.17). The 1/2 *sar* rectangle is complete, except for a small rectangular piece.

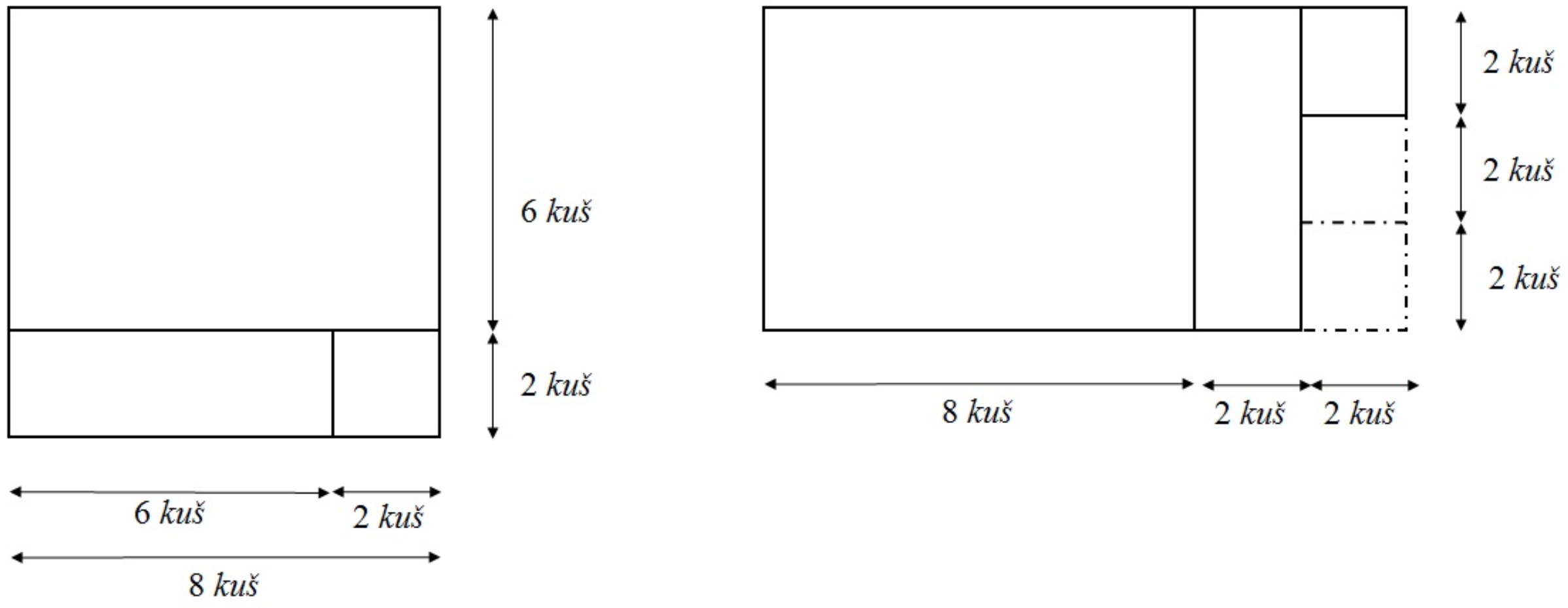

**Fig. 9.17** Surface of the 8 *kuš*-side square

The elementary squares and rectangles are that measuring 1/2 *sar* (sides 12 *kuš* and 6 *kuš*), and the missing rectangular piece (sides 4 *kuš* and 2 *kuš*, that is, two 2 *kuš*-side squares); these elementary surfaces are represented below (Fig. 9.18):

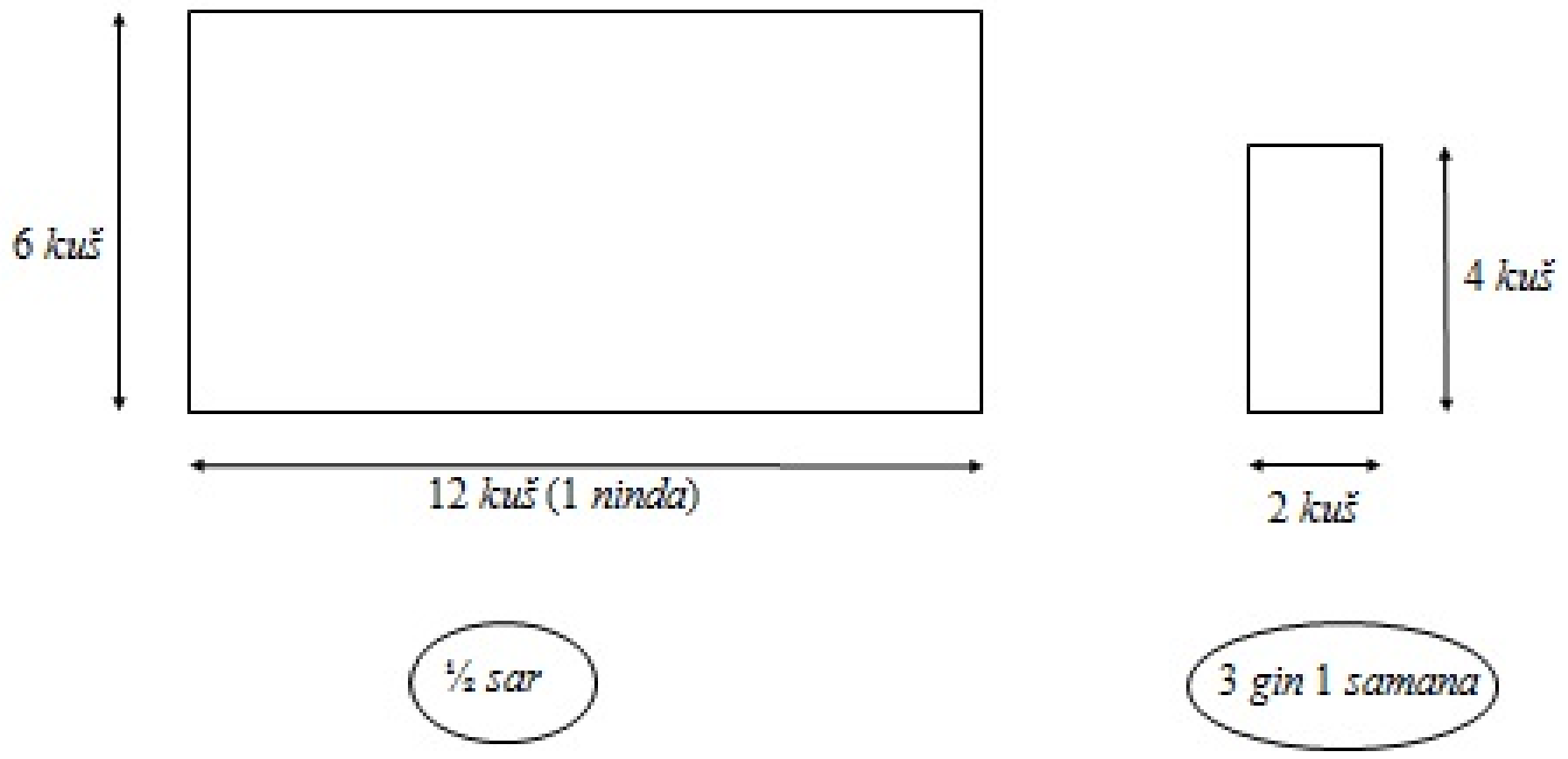

**Fig. 9.18** Elementary surfaces (continuation)

The surface noted '1/2 *sar* – (3 *gin* 1 *samana*)' represents exactly the large piece from which the small piece is cut. Here again, the pieces which appear in the subtractive notations represent elementary surfaces.

If this interpretation is accepted, the quantification of small surfaces relies on diverse arrangements of elementary surfaces formed from the measurement units. This process is more complex than in the case for those used in Texts 1-3 because the elements of surfaces are not only pasted, but also cut, leading to subtractive notations for the surfaces.

The process, as described above, is based on the preliminary knowledge of a repertoire of elementary surfaces, for example, the result provided by the first clause of the text (the surface of the 1 *kuš*-square is 1 *samana* 15 *še*). Moreover, the second clause should be obtained from the first one in a very simple way: if the surface of a 1 *kuš*-square is 1 *samana* 15 *še*, then the surface of a 2 *kuš*-square is four times more, that is, 1 *gin* 2 *samana* (4 *samana* 60 *še*). But the text gives the result in the unexpected subtractive form ‘2 *gin minus* 1 *samana*’.

The previous observations (the omission of the 1 *ninda* entry, corresponding to the obvious surface 1 *sar*; the use of previous non-evident results, such as the surface of a 1 *kuš*-square is 1 *samana* 15 *še*; the expression of the surfaces in unexpected subtractive forms) suggest that this text was not conceived as a list of results to be used for evaluating lands, as Texts 1-3, but as a list of mathematical problems. The statements of these problems are not explicitly formulated, they can hardly have been a simple question such as ‘if the side of a square is this, what is the surface?’, because the answer of the second problem should be ‘1 *gin* 2 *samana*’ and not ‘2 *gin minus* 1 *samana*’. Thus, the statements of the problems refer to a complex process of cutting and pasting pieces of surfaces which left no trace, except the subtractive notations.

The tentative reconstruction of procedures for evaluating surfaces above suggests that some high-ranked land-surveyors (or scholars linked to these land-surveyors) developed, in Early Dynastic IIIa and b periods, a sophisticated technology of cutting and pasting pieces of surfaces. This knowledge may have stimulated a long mathematical tradition, still alive half a millennium later in the Old Babylonian period, when geometrical methods were developed for solving quadratic problems (Høyrup 2002).

Another solution to the problem of computing small surfaces, completely different from the metrological manipulations described above, emerges from the examination of Text 5.

### *9.4.3 Text 5 (CUNES 50-08-001)*

Tablet CUNES 50-08-001, published by Friberg (2007: appendix 7) is thought to come from Zabalam and is dated to the Early Dynastic IIIb (*ibid*: 419). The document is a large multi-column tablet with seven columns on the obverse, and three columns followed by a large blank space on the reverse. The text ends with quite a long subscript (rev. *iii* 3-9), the meaning of which is not clear. Friberg (*ibid.*: 422-425) provided a sketch of transliteration and detailed explanations; a complete transliteration and translation can be found in Appendix 9.C.

CUNES 50-08-001 contains a set of five sub-tables. Each sub-table is a list of clauses providing the sides of squares and the corresponding surfaces listed in increasing size order. The sides of the squares cover a vast range, from 1 *šu-bad* (ca. 25 cm) to 36 000 *ninda* (ca. 216 km). Thus, we can consider that this table does not deal with fields, like Texts 1-3, but with abstract squares, the dimensions of which can be as large or small as the numerical and metrological systems allowed. Furthermore, the available metrological system for lengths seem to have been extended with the use of unusually large System S numbers (*šaru* and *šar*, see Fig. 9.2), and a tiny sub-division of *ninda* (*šu-bad*).

A unique feature of CUNES 50-08-001 is the metrological system used for small surfaces. As shown by Friberg, the scribes used an extended system of sexagesimal sub-units of the *sar* unit, which is not attested elsewhere. The metrology used in Text 5 is represented by the following diagram (Table 9.6) and will be explained in more details below.

| | | |
|---|---|---|
| Units of length: | | ninda-DU ←×4— nig$_2$-kas$_7$ ←×1 ½— kuš$_3$-numun ←×2— giš-bad ←×2—šu-bad |
| | | ↕ |
| Units of surface: | GAN$_2$ ←×100— sar ←×60— gin$_2$ ←×60— gin$_2$-bi ←×60— gin$_2$-ba-gin$_2$ | |

**Table 9.6** Measurement units used in CUNES 50-08-001

### *9.4.4 Organisation of Text 5*

The five sub-tables of Text 5, numbered as A, B, C, D, E by Friberg (2007: 419), are represented in Fig. 9.19 below.

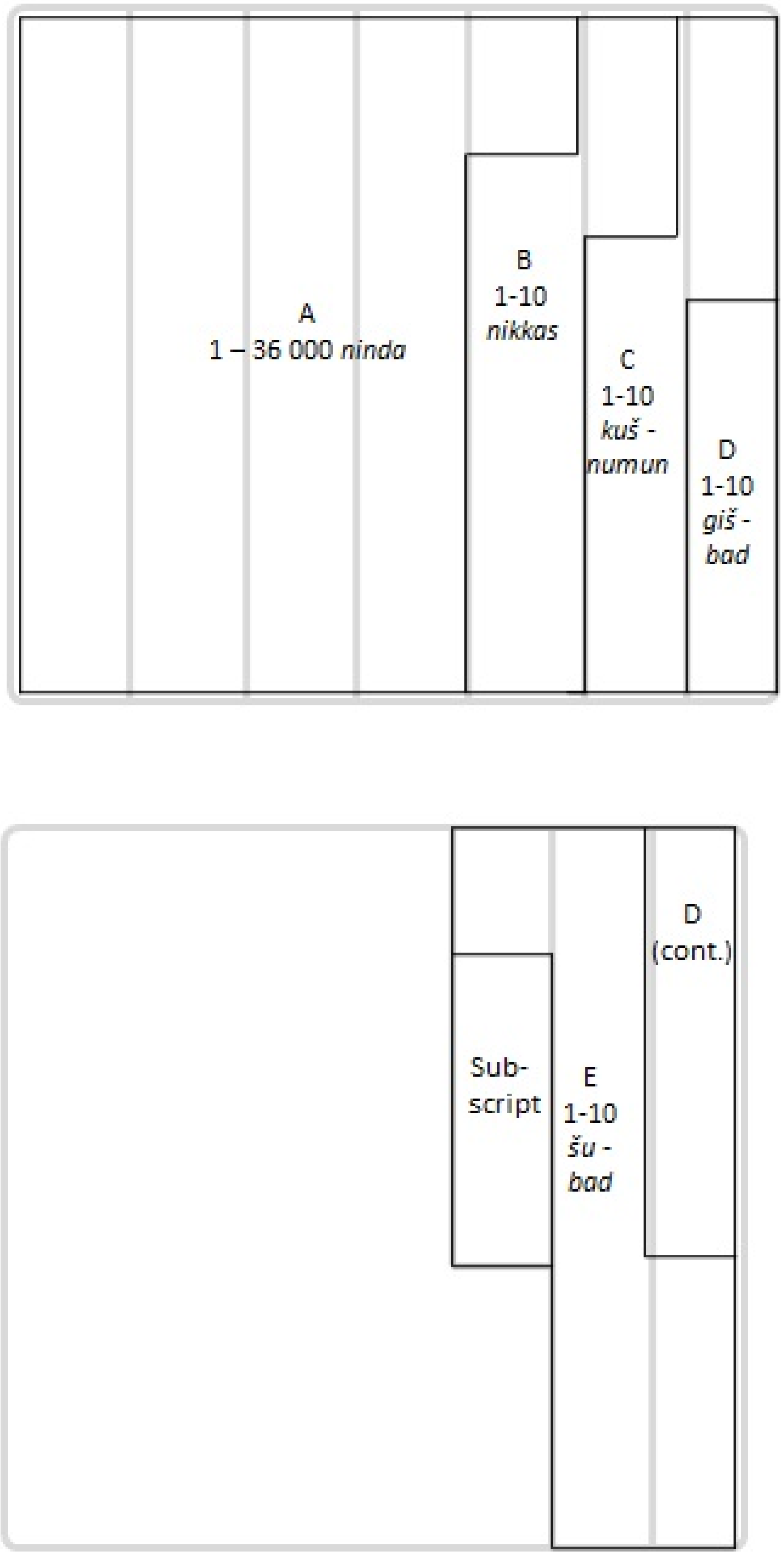

**Fig. 9.19** Organisation of Text 5

As shown by Friberg (*ibid*), sub-table A gives the surfaces of *N ninda*-side squares (*N* integer from 1 to 36 000). Then, sub-tables B, C, D and E give the surfaces of squares for each sub-unit of the *ninda*. The content of the sub-tables can be summarised as follows:

- 1 to 36 000 *ninda* in sub-table A
- 1 to 10 *nikkas* in sub-table B
- 1 to 10 *kuš-numun* in sub-table C
- 1 to 10 *giš-bad* in sub-table D
- 1 to 10 *šu-bad* in sub-table E.

The five sub-tables are listed in decreasing size order, from lengths in *ninda* to lengths in *šu-bad*, while the lengths inside each sub-table are listed in increasing size order. The uniform range, from 1 to 10, adopted for numerical values in all of the tables B-D doesn't reflect the diversity of factors in the metrological system of lengths (4, 1 1/2, 2, and 2, see Table 9.6). As a consequence, Text 5 contains many redundant items. For example, the lengths 2 *nikkas* (in sub-table B, obv. *v* 6), 3 *kuš-numun* (in sub-table C, obv. *vi* 11) and 6 *giš-bad* (in sub-table D, obv. *vii* 18) are equal, thus, the corresponding surfaces are the same (1 *sar* 15 *gin*), and this surface appears three times in the text (obv. *v* 7; obv.*vi*, 12; and rev. *i*, 1). Table 9.7 below lists the redundancies (each row contains the equal lengths appearing in the different sub-tables).

| A | B | C | D | E |
|---|---|---|---|---|
| *ninda* ←4— | *nikkas* ←1 1/2— | *kuš-numun* ←2— | *giš-bad* ←2— | *šu-bad* |
| | | | 1 *giš-bad* | 2 *šu-bad* |
| | | 1 *kuš-numun* | 2 *giš-bad* | 4 *šu-bad* |
| | 1 *nikkas* | | 3 *giš-bad* | 6 *šu-bad* |
| | | 2 *kuš-numun* | 4 *giš-bad* | 8 *šu-bad* |
| | | | 5 *giš-bad* | 10 *šu-bad* |
| | 2 *nikkas* | 3 *kuš-numun* | 6 *giš-bad* | |
| | | 4 *kuš-numun* | 8 *giš-bad* | |
| | 3 *nikkas* | | 9 *giš-bad* | |
| | | 5 *kuš-numun* | 10 *giš-bad* | |
| 1 *ninda* | 4 *nikkas* | 6 *kuš-numun* | | |
| | 6 *nikkas* | 9 *kuš-numun* | | |
| 2 *ninda* | 8 *nikkas* | | | |

**Table 9.7** Redundant information in Text 5

These features (the uniform range of numerical values and the resulting repetitions of the same information) suggest that the main purpose of Text 5 was not to provide information, but rather to exhibit mathematical ideas about surfaces. These ideas may emerge from the reconstruction of the methods used by the ancient scribes for computing surfaces.

### *9.4.5 Computation*

#### Sub-table A

Was the method of computation of large surfaces (sub-table A) the same as in Texts 1-3 analysed above? As most of the data and notations are the same, one may be tempted to answer in the affirmative. However, not only the format (tabular in Texts 1-3, list of clauses in Text 5), but also two other important features contrast Texts 1-3 and sub-table A of Text 5: in Texts 1-3, the smallest length is 5 *ninda*, and the *sar* unit is not used; in sub-table A of Text 5, the smallest length is 1 *ninda* and the *sar* unit is used. Thus, the bridge between the units of length and the units of surface seems to have moved from 10 *ninda* (with corresponding unit of surface 1 *iku GAN*) in Texts 1-3, to 1 *ninda* (with corresponding unit of surface 1 *sar*) in sub-table A of Text 5. Thus, the structure of the metrological system adopted in table A (see Table 9.8) differs from the one adopted in Texts 1-3.

| | | |
|---|---|---|
| Unit of length: | | ninda-DU |
| | | ↕ |
| Unit of surface: | GAN$_2$ ←×100− | sar |

**Table 9.8** Measurement units used in Text 5 (CUNES 50-08-001), sub-table A

As the articulation between length and surface measures is not the same in ED IIIa and ED IIIb tables, the calculation method may have differed. The hypothesis we have rejected for tables 1-3 may be more probable for sub-table A, Text 5: for each field, the number of *ninda* measuring the side (*N*) may have been multiplied by itself. The resulting product (*N*×*N*) represents a number of *sar*, and would have to be converted into System G. This arithmetical procedure is consistent with the arithmetical approaches to the surfaces of small squares which emerge from the following analysis of sub-tables B-E.

**Small squares**

Friberg (2007: 422-4) reconstructed the method of computation for small surfaces (sub-tables B, C, D, E) supposing the use of a geometrical variant of the binomial rule. The following attempt is not fundamentally different from Friberg's suggestions. However, it tries to rely more closely on the specific features of the text. We just observed that the entries are organised in a systematic way, by considering each sub-unit of the *ninda*, and entering 1, 2, …, 10 of these sub-units. We can also observe that the notations of surfaces do not use any subtractive notation. This may imply that the surfaces were not obtained by geometrical cut-and-paste procedures similar to the ones described for Text 4. The most striking feature of this text is the systematic use of sexagesimal sub-divisions of the *sar*, namely a sixtieth of a *sar*, a sixtieth of a sixtieth of a *sar*, and a sixtieth of a sixtieth of a sixtieth of a *sar*, translated by Englund (CDLI: P274845) as 'shekel', 'small-shekels', and 'mini-shekels'.

Not only is such an extended system of sexagesimal sub-divisions of measurement units quite unusual,[32] but also the notations associated to these sexagesimal *sar* sub-divisions are puzzling, and it was certainly a challenge to elucidate the meaning of tables B-D (see correspondence with P. Damerow in CDLI: P274845). These subdivisions appear, for example, in the surface corresponding to the first entry in table B, that is, the surface of a 1 *nikkas*-side square (Fig. 9.20).

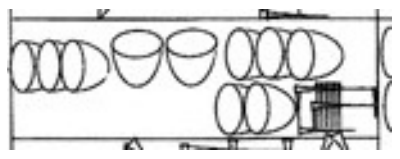

3(aš) 2/3 5(aš) gin$_2$

**Fig. 9.20** Text 5 (CUNES 50-08-001), table B, obv. *v* 5, copy Friberg 2007: 420-421

The computation of the surface of a 1 *nikkas*-side square gives 3 2/3 sixtieths of a *sar* and 5 sixtieths of a sixtieth of a *sar*. This is actually the surface given by the text, but the notation is quite defective: the numerical values are noted, but the measurement units are partially implicit, which gives the notation an almost positional appearance.

In other sections, the names of the measurement units appear, but in an unusual disposition. For example, the surface corresponding to 1 *giš-bad*, the first entry in table D is as follows (Fig. 9.21).

32 Friberg (2007: 419) mentions a handful of Old Akkadian texts exhibiting such sexagesimal sub-units of the *sar*. However, the sexagesimal system used in these texts is not so developed as in our Early Dynastic text.

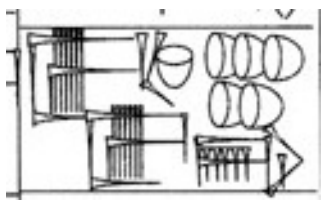

gin$_2$ 1/3 5 gin$_2$-bi

**Fig 9.21** Text 5 (CUNES 50-08-001), table D, obv. *vii* 10, copy Friberg 2007: 420-421

The computation of the surface of a 1 *giš-bad*-side square gives 1/3 of a sixtieth of a *sar* and 5 sixtieth of a sixtieth of a *sar*. As the sign gin$_2$ was widely used with the meaning of a sixtieth in subsequent Sumerian metrology, we can deduce that the first sign represents a sixtieth of a *sar*, and thus, the final signs 'gin$_2$-bi' represent a sixtieth of a sixtieth of a *sar*. Thus, the notation above (Fig. 9.21) can be understood 1/3 *gin* 5 *gin-bi*. The further we advance in the text, the stranger the notations. For example, the surface corresponding to 1 *šu-bad*, the first entry in table E is as follows (Fig. 9.22).

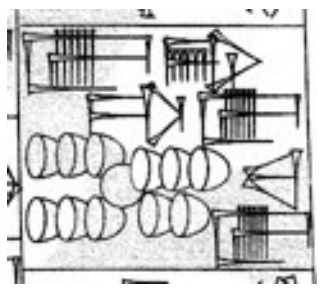

gin$_2$-bi-ta 6(aš) 1(u) 5(aš) gin$_2$-ba-gin$_2$

**Fig. 9.22** Text 5 (CUNES 50-08-001), table E, rev. *i* 11, copy Friberg 2007: 420-421

The computation of the surface of a 1 *šu-bad*-side square gives 6 sixtieth of a sixtieth of a *sar* and 15 sixtieth of a sixtieth of a sixtieth of a *sar*. We can understand the sequence 'gin$_2$-ba-gin$_2$' as representing a sixtieth of a sixtieth of a sixtieth of a *sar* (1 mini-shekel), and the notation of Fig. 9.22 above as 6 *gin-bi* 15 *gin-ba-gin* (6 small-shekels 15 mini-shekels). This is Friberg's interpretation, but Englund, according to his translation in CDLI ('of its shekel 6 small-shekels, 15 [mini-shekels] its small-shekel'), seems to understand the syntax of the notation in a slightly different way. Despite the difficulty in understanding this syntax, I adopt the representations of the sexagesimal sub-units of the *sar* suggested by Friberg:

1 *gin* = 1/60 *sar*
1 *gin-bi* = 1/60 *gin*
1 *gin-ba-gin* = 1/60 *gin-bi*

Let us come back to the computations. For each table the computation of the surface of the first square, that is, the square whose side is 1 sub-unit of the *ninda* (respectively 1 *nikkas*, 1 *kuš-numun*, 1*giš-bad* and 1 *šu-bad*) is examined. Then, the other surfaces of each sub-table are considered.

**Sub-table B**

The first entry is a 1 *nikkas*-side square. As 1 *nikkas* is a quarter of 1 *ninda*, the surface of a 1 *nikkas*-side square is a sixteenth of the surface of a 1 *ninda*-side square, that is, a sixteenth of 1 *sar*. Thus, the starting item in sub-table B is based on the preliminary knowledge that:

1/16 *sar* = 3 2/3 *gin* 5 *gin-bi* (1)

The following two items in sub-table B can have been obtained by multiplying this result, 3 2/3 *gin* 5 *gin-bi*, by 4 and 9 respectively. A simple additive procedure may have led to the results 15 *gin* and 1/2 *sar* 3 2/3 *gin* 5 *gin-bi* respectively (see Table 9.9). The following item provides the surface of a 4 *nikkas*-side square, that is, a 1 *ninda*-side square, which is known to be 1 *sar*. The 8 *nikkas*-side square (=2 *ninda*-side square) is four times the 1 *ninda*-side square, that is, 4 *sar*. The other item may have been obtained by additive procedures similar to those used for the 2 and 3 *nikkas*-side squares.

| **Equal side (sa$_2$)** | **Surface noted on the tablet** | **Calculated surface** |
|---|---|---|
| 1 *nikkas* | 3 2/3 5 gin$_2$ | 3 2/3 *gin* 5 *gin-bi*<br>(= 3 *gin* 45 *gin-bi*) |
| 2 *nikkas* | 15 gin$_2$ | 15 *gin* |
| 3 *nikkas* | 1/2 sar 3 2/3 5 gin$_2$ | 1/2 *sar* 3 2/3 *gin* 5 *gin-bi*<br>(= 33 *gin* 45 *gin-bi*) |
| 4 *nikkas* (=1 *ninda*) | 1 sar | 1 sar |
| 5 *nikkas* | 1 1/2 sar 3 2/3 5 gin$_2$ | 1 1/2 *sar* 3 2/3 *gin* 5 *gin-bi*<br>(=1 *sar* 33 *gin* 45 *gin-bi*) |
| 6 *nikkas* | 2 15 gin$_2$ | 2 *sar* 15 *gin* |
| 7 *nikkas* | 3 sar 3 2/3 5 gin$_2$ | 3 *sar* 3 2/3 *gin* 5 *gin-bi*<br>(= 3 *sar* 3 *gin* 45 *gin-bi*) |
| 8 *nikkas* (=2 *ninda*) | 4 sar | 4 *sar* |
| 9 *nikkas* | 5 sar 3 2/3 5 gin$_2$ | 5 *sar* 3 2/3 *gin* 5 *gin-bi*<br>(= 5 *sar* 3 *gin* 45 *gin-bi*) |
| 10 *nikkas* | 6 sar 15 gin$_2$ | 6 *sar* 15 *gin* |

**Table 9.9** Calculation of surfaces in sub-table B

As we see, only the first item, based on relation (1) cannot be deduced from previous items of the text by additive procedures, and thus requires preliminary knowledge.

**Sub-table C**

The first entry is a 1 *kuš-numun*-side square. As 1 *kuš-numun* is a sixth of 1 *ninda*, the surface of a 1 *kuš-numun*-side square is a thirty-sixth of the surface of a 1 *ninda*-side square, that is, a thirty-sixth of 1 *sar*. Thus, the starting item of sub-table C is based on the preliminary knowledge that:

1/36 *sar* = 1 2/3 *gin* (2)

The surface of the 6-*kuš-numun*-side square (= 1 *ninda*-side square) measures 1 *sar* evidently. The surface of the 3-*kuš-numun*-side square (= 1/2 *ninda*-side square) measures 1/4 *sar* evidently, which is noted as 15 *gin* in the text. The other surfaces in sub-table C are obtained by multiplying 1 2/3 *gin* by 4, 16, …, 100 respectively, possibly with additive procedures.

**Sub-table D**

The first entry is a 1 *giš-bad*-side square. As 1 *giš-bad* is half of 1 *kuš-numun*, the surface of a 1 *giš-bad* -side square is a quarter of the surface of a 1 *kuš-numun* -side square, that is, quarter of 1 2/3 *gin*. Here, a new sexagesimal sub-division appears: 1 *gin-bi* is 1/60 *gin*. Thus, a quarter of 1 2/3 *gin* is a quarter of 100 *gin-bi*, that is, 25 *gin-bi* or, equivalently, 1/3 *gin* 5 *gin-bi*. Thus, the starting item in sub-table D is based on the preliminary knowledge that:

1/4 (1/36 *sar*) = 1/3 *gin* 5 *gin-bi* (3)

The surface of a 6 *giš-bad*-side square (= 1/2 *ninda*-side square) measures 1/4 *sar* evidently, that is, 15 *gin*. The other surfaces in sub-table C are obtained by multiplying 1/3 *gin* 5 *gin-bi* by resp. 4, 9, 16, …, 100, possibly with additive procedures.

**Sub-table E**

The first square is 1 *šu-bad*-side. As 1 *šu-bad* is a half of 1 *giš-bad*, the surface of a 1 *šu-bad*-side square is a quarter of the surface of a 1 *giš-bad*-side square, that is, a quarter of 25 *gin-bi*. Here, a new sexagesimal sub-division appears: 1 *gin-ba-gin* is 1/60 *gin-bi*. Thus, a quarter of 25 *gin-bi* is a quarter of 24 *gin-bi* 60 *gin-ba-gin*, that is, 6 *gin-bi* 15 *gin-ba-gin*. Therefore, the starting item in sub-table E is based on the preliminary knowledge that:

1/4 (25 *gin-bi*) = *6 gin-bi 15 gin-ba-gin* (4)

The other surfaces in sub-table E are obtained by multiplying = 6 *gin-bi* 15 *gin-ba-gin* by resp. 4, 9, 16, …, 100, possibly with additive procedures.

The construction of sub-tables B-E seems to be based on the preliminary knowledge summarised in Table 9.10.

| Sub-table | Square | Preliminary knowledge | | |
|---|---|---|---|---|
| B | 1 *nikkas* □ 1 *nikkas* | 1/16 *sar* | = 3 2/3 *gin* 5 *gin-bi* | (1) |
| C | 1 *kuš-numun* □ 1 *kuš-numun* | 1/36 *sar* | = 1 2/3 *gin* | (2) |
| D | 1 *giš-bad* □ 1 *giš-bad* | ¼ (1/36 *sar*) | = 1/3 *gin* 5 *gin-bi* (=25 *gin-bi*) | (3) |
| E | 1 *šu-bad* □ 1 *šu-bad* | ¼ (25 *gin-bi*) | = 6 *gin-bi* 15 *gin-ba-gin* | (4) |

**Table 9.10** Preliminary knowledge in the construction of sub-tables B-E

This preliminary knowledge reflects a mastery of equivalencies between fractions of *sar*, and sexagesimal subdivisions of *sar*. Thus, we can hypothesise that a complete set of equivalencies such as Table 9.11 below was known by the author of Text 5:

| | | | |
|---|---|---|---|
| 1 *sar* | = | 60 *gin* | |
| 1/2 *sar* | = | 30 *gin* | |
| 1/3 *sar* | = | 20 *gin* | |
| 1/4 *sar* | = | 15 *gin* | |
| 1/5 *sar* | = | 12 *gin* | |
| 1/6 *sar* | = | 10 *gin* | |
| 1/8 *sar* | = | 7 1/2 *gin* | |
| 1/9 *sar* | = | 6 2/3 *gin* | |
| 1/10 *sar* | = | 6 *gin* | |
| 1/12 *sar* | = | 5 *gin* | |
| 1/15 *sar* | = | 4 *gin* | |
| 1/16 *sar* | = | 3 2/3 *gin* 5 *gin-bi* | (1) |
| […] | | | |
| 1/36 *sar* | = | 1 2/3 *gin* | (2) |
| […] | | | |
| ¼ (1/36 *sar*) | = | 1/3 *gin* 5 *gin-bi* | (3) |
| […] | | | |
| 1/16 (1/36 *sar*) | = | 6 *gin-bi* 15 *gin-ba-gin* | (4) |

**Table 9.11** Fractions of a *sar* and sexagesimal sub-divisions

Moreover, the generation of each sub-table from the starting item seems to involve some ability to multiply a sexagesimal expression by an integer. This skill may have been based on additive procedures, as suggested above, or on the knowledge of some multiplication tables, or on the manipulation of some sort of device.

The process of calculation tentatively reconstructed above may explain the structure of the text. First, table A of large fields (sides of 1 *ninda* and beyond) is computed, perhaps by arithmetical procedures (calculating the arithmetical square of a number). Then, the 1 *nikkas*-side square is calculated, possibly using a table of correspondences between fractions of a *sar* and sexagesimal sub-divisions of a *sar* such as Table 9.11. Table B for 2 *nikkas*, 3 *nikkas*, …, 10 *nikkas* is obtain by multiplying the 1 *nikkas*-side surface by resp. 4, 9, …, 100, possibly using multiplication tables or a device. The same process, which introduces new sexagesimal sub-divisions of the units of surface when necessary, produces tables C, D, and E. The unexpected subdivisions of the *ninda* unit adopted in this text, and not attested elsewhere, may have been created for the exploration of the newly discovered sexagesimal methods of computation.

Text 5 turns out to be a systematic exploration of sexagesimal computation. In this reconstruction, the key tool would be a table of correspondence similar to Table 9.11, which was later to become, after the invention of sexagesimal place value notation, a reciprocal table.[33]

## 9. 5 Conclusion

Different approaches to the quantification of surfaces, involving different sorts of multiplication, emerge from this overview of Early Dynastic tables. Table 9.12 summarize the way in which these different approaches correlate with the diverse features of the tables.

| **Tables** | **Size of surf.** | **Layout** | **Notations of surfaces** | **Metrology** |
|---|---|---|---|---|
| 1 | Large | Tabular | Simple | 10 ninda-DU<br>↕<br>GAN$_2$ |
| 2, obv. | Large | Tabular | Simple | 10 ninda-DU<br>↕<br>GAN$_2$ |
| 3, A | Large | Tabular | Simple | 10 ninda-DU<br>↕<br>GAN$_2$ |
| 3, B | Large | Tabular | Simple | 10 ninda-DU ←10— ninda-DU ←4− ur$_2$ hal-la ←1 1/2 − kuš$_3$-n.<br>↕<br>GAN$_2$ |
| 4 | Small | List | Compound; subtractive notations | (ninda ←2−) gi ←6− kuš$_3$<br>↕<br>sar ←60− gin$_2$ ←3− sa$_{10}$-ma-na ←60− še |
| 5, A | Large | List | Simple | ninda ←4— nig$_2$-k. ←1 1/2— kuš$_3$-n. ←2— giš-b. ←2— šu-b.<br>↕<br>GAN$_2$ ←100− sar ←60— gin$_2$ ←60— gin$_2$-bi ←60— gin$_2$-ba-gin$_2$ |
| 5, B, C, D, E | Small | List | Compound; sexagesimal sub-divisions | |

**Table 9.12** Features of the Early Dynastic tables of surfaces

Throughout this analysis, I contrasted large lands with small squares and rectangles. The reason for this demarcation becomes clear when considering the metrologies represented in the last column of Table 9.12. The surfaces of large lands, which appear in tables 1, 2, 3 A-B, are greater than 1/4 *iku GAN* (= 25 *sar*). The pivot of the relationship between lengths and surface measurement values is the 10 ninda-side square (surface 1 *iku GAN*). The procedure for the calculation of surfaces seems to have been based on metrological manipulation.

The surfaces of small squares and rectangles, which appear in Text 4, are all under 1 or 2 *sar*. The pivot of the relationship between lengths and surface measurement values is the 1 *ninda*-side square (surface 1 *sar*). This means that the units used for sides and surfaces are downstream of the bridge between length and units of surface—as in the other representations of metrological systems in this chapter, in Table 9.12 these 'bridges' are represented by double vertical arrows (↕). Thus, the evaluation of the surfaces involves small sub-units of the *sar*, and it is impossible to calculate the surface using an arithmetical multiplication acting on integers. But it seems that in this text, there is no hint of the use of calculation with fractions. The facts that the notation of the surfaces use subtractive notation, and that the measurement values are compound may show that the computations are based on cutting and pasting elementary pieces of surfaces.

[33] The link between the sexagesimal sub-divisions of *sar* in CUNES 50-08-001 and the subsequent invention of the sexagesimal place value notation was highlighted by Friberg (2007: 426).

Text 5 testifies a much more unified approach to surfaces. All types of fields are presented in the same tablet, whether large or small. This highly structured set of tables shows a shift from geometrical to arithmetical procedures to evaluate the surfaces. Sub-table A covers a very large range of fields (sides from 1 to 36 000 *ninda*) and may have been computed using arithmetical multiplication. In sub-tables B-E, the notations of surfaces do not use subtractive notation, but are compound. The unit of surface '*sar*' is divided into sexagesimal sub-units, which are not attested elsewhere and appear to have been invented to help in handling fractions of *sar*. The computation of surfaces seems to be based on the knowledge of some elementary results of sexagesimal computation (sexagesimal equivalents of fractions of *sar*), that is, a kind of non-positional version of a reciprocal table.[34]

The different procedures for quantifying surfaces can be summarized as follows:

1) The quantification of surfaces in Texts 1-4 seems to rely on diverse arrangements of elementary surfaces formed from the measurement units. In this case, the operations are geometrical and act on units of length, and these operations produce elementary surfaces such as those represented in Fig. 9.11 and 9.12. In the case of small squares (Text 4), these metrological manipulations explain why the resulting surfaces are noted with subtractive notations or compound measurement values.
2) The quantification of large surfaces in sub-table A in Text 5 seems to be based on arithmetical multiplication on non-positional sexagesimal numbers.
3) The quantification of small surfaces in sub-tables B-E in Text 5 is clearly based on sexagesimal computation. The operations include taking fractions of fractions of units (e.g. a quarter of the quarter of a unit is the sixteenth of this units), and then using correspondences between fractions and sexagesimal sub-divisions (such as those provided in Table 9.11). The resulting surface is noted with compound measurement values.

The function of the tables in Texts 1-3 seems to be different from the function of the tables in Texts 4-5. Texts 1-3, in a tabular format, deal with large lands, and appear to be mathematical elaborations inspired by surveying practices. By contrast, some of the surfaces in Texts 4-5 have unrealistic sizes (much smaller, or much larger than real plots and lands), and the layout of the texts is a list of clauses. Text 4 seems to denote a reflexion, expressed by lists of solved problems, on the notion of surface. This notion, applied to small squares, emerges as a theoretical extension of a common notion of surface applied by land-surveyors to fields and lands. This extension turned out to raise a difficult mathematical problem due to the fact that there is only one bridge between units of length and units of surface. Text 5 can be interpreted as an attempt to both unify the different approaches to surface and to solve the problem of small surfaces, by mean of arithmetical tools based on sexagesimal computation.

---

[34] Friberg (2019: Sect. 10) published a non-positional table of reciprocasl, SM 2685, dated to the Ur III or early OB period (end of the third millennium or very beginning of the second). Although much later than the texts considered here, this table may echo some tables of equivalencies between fractions and sexagesimal sub-divisions possibly used formerly in Early-Dynastic or Sargonic periods.

## Appendix 9.A: System G and the Sign GAN$_2$

The notation of surface measurements in early third millennium documentation, and the distinction between a numerical value and a measurement unit raised by this notation, were examined by Powell (1972, 1973), who formulated the problem as follows:

> The Sumerian measures of area conform very closely to the system of numeration both in notation and in actual structure. … An important point to bear in mind in trying to understand the System of area measures is that *metrological notation is made by the same graphic procedure used to make numerical notation.* (his emphasis)

According to Powell, the signs used to quantify surfaces are both numerical and metrological. In other words, Powell assumes that there is no dissociation between numerical value and measurement unit in the notation of measurement values of surfaces. However, it seems to me that, while this assumption is relevant for the earliest texts using System G, the situation is much more ambiguous and diverse in the context of the Early Dynastic tables discussed here.

In the quotation above, Powell refers to documentation dealing only with surfaces of large fields using System G and dated to earlier periods than the tables examined here. By contrast, the notation of surfaces with *sar*, a smaller unit of surface (ca. 36 m²), which perhaps appeared later, fits perfectly with the general arrangement 'numerical value + measurement unit' drawn above, as shown by Fig. 9.23.

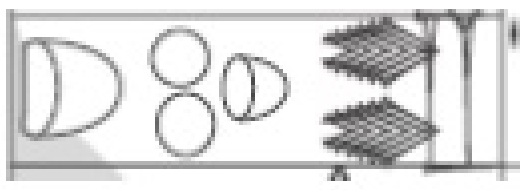

1(geš$_2$) 2(u) 1(aš) sar
(1×60 + 2×10 + 1 *sar*)

**Fig. 9.23** Text 5 (CUNES 50-8-1) obv. i 18, extracted from Friberg (2007: 420 Fig. A7.1)

In this example, the signs '1(geš$_2$) 2(u) 1(aš)' make up a numerical value according to System S, and the sign 'sar' is a metrogram which denotes a measurement unit belonging to the metrological system of surfaces. The dissociation between numerical and metrological components in the notation of surfaces is clearly accomplished in the notations using *sar* (and sub-divisions of *sar*).

In fact, Powell's analysis quoted above refers only to surface measurements associated with System G and the sign GAN$_2$. The discussion can be focused on the question: how should a notation such as the following be analysed (Fig. 9.24)? [35]

[35] Somehow the answer is given by the translation that I have adopted in Fig. 9.23-24, and that I shall now attempt to justify. Actually, my transliterations and translations of surfaces involving the sign GAN$_2$ are quite close to the conventions used by Neugebauer (see for example Neugebauer 1935: 91). However, the conventions used in modern publications are quite different, and this situation deserves some clarification.

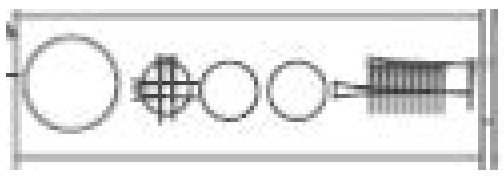

1(šar$_2$) 1(bur'u) 2(bur$_3$) GAN$_2$
(1080 + 180 + 2×18 GAN$_2$)
**Fig. 9.24** Text 5 (CUNES 50-8-1) obv. ii 18, taken from Friberg (2007: 420 Fig. A7.1)

The first problem concerns the sign GAN$_2$. Powell shows in the two publications cited above (1972, 1973) that the function, pronunciation and meaning of the sign GAN$_2$ in third millennium documentation are far from clear and, actually, are not uniform. He argues that in certain contexts, the sign GAN$_2$ should be read as 'aša' and means 'field' or 'area' (Powell 1973: 182). In other contexts the sign GAN$_2$ should be read as 'gana' and means 'land, ground, soil or the like' (Powell 1973: 183). And finally, in notations of surface measurements, the sign GAN$_2$ is a simple semantic indicator not to be read at all (Powell 1973: 182), such as in Early Dynastic IIIb Lagaš texts (Nissen, Damerow and Englund 1993: 64).

In our Early Dynastic tables, does the sign GAN$_2$ act as a measurement unit, or as the quantified quantity ('field'), or as a simple semantic indicator, that is, a graphic mark with no verbal counterpart? Or is the sign ambiguous, incorporating some or all of these elements (to greater and lesser extents)?

The presence of the sign GAN$_2$ in the heading of tabular tables (see Text 2) in a similar way as the unit of length *ninda*, advocates for the first hypothesis. Moreover, in our texts, the sign GAN$_2$ always appears when a unit of surface corresponding to 100 *sar* is expected. For example, GAN$_2$ never follows a simple measurement value expressed in *sar* (something like 10 *sar* GAN$_2$ does not exist). Thus, GAN$_2$ does not denote a field in general, but a special field whose surface measures 100 *sar*, or, in texts where the unit *sar* does not appear, a 10 *ninda*-side square. From a functional point of view, the sign GAN$_2$ acts as a unit of surface, even if other functions suggested above are also present.

The second problem concerns the status of notations such as '1(šar$_2$) 1(bur'u) 2(bur$_3$)' in Fig. 9.24. In the two examples above (Fig. 9.23 and 9.24), the structure of the measurement values is emphasized through writing techniques: curve signs, made by scribes with the rounded end of a reed, are followed by a cuneiform sign, made with the sharp end of a reed. Thus, one is tempted to recognize the same structure in the second example as in the first. In other words, in the same way as the signs '1(geš$_2$) 2(u) 1(aš)' form a numerical value and are followed by a cuneiform sign representing 1 *sar* (Fig. 9.23), the signs '1(šar$_2$) 1(bur'u) 2(bur$_3$)' form a numerical value and are followed by a cuneiform sign, GAN$_2$, representing 100 *sar* (Fig. 9.24). In this latter notation (Fig. 9.24), one can recognize the two components 'numerical value + measurement unit', where the numerical values belong to the so-called 'System G' (see Fig. 9.4 in Sect. 9.1), and the measuring unit is GAN$_2$. This interpretation leads to consider that System G is a kind of numerical system. In this line, one can observe that the syntax of System G is the same as System S. Indeed, each sign is noted as many times as necessary, which is characteristic of additive numerical systems. The comparison of the signs adopted in Systems S and G (see Fig. 9.2,9.4 and Sect. 9.1) shows that the graphical repertory is partly identical, which reinforces the temptation to consider System G as a numerical system. Under these assumptions, numerical values would be visually perceptible in our texts by the fact they all of them, whether they belong to System S or G, are noted as curve signs, while measurement units and the names of the quantified quantities are noted as cuneiform signs. In Texts 1 and 2, these numerical signs appear in all the items of the different tabular columns, while the measurement units and the quantified quantity only appear in headings.

The semantic of disposition in Early Dynastic tables tends to assign a numerical function to signs in System G, and a function close to a measurement unit to the sign $GAN_2$. This means that Early Dynastic tables show a clear attempt to separate the numerical components from the metrological components in the notations of measurement values. In this shift, the function of the sign $GAN_2$ seems to have changed, perhaps in connection with the introduction of the smaller *sar* unit. This attempt may have been local, limited to some erudite scribal milieus, and have not affected other contemporary practices, for example in administration (see for example the field texts from Lagaš analysed by Camille Lecompte in this volume). But this attempt certainly had an impact on subsequent mathematical writings produced in scribal schools, as the separation between the numerical and metrological components is clearly accomplished for surface measures in Old Babylonian mathematical texts.

## Appendix 9.B: Chronology

Estimated dates adopted by the CDLI, according to the middle chronology.

| **Periods** | **Dates BCE** |
|---|---|
| Uruk IV | ca. 3350-3200 |
| Uruk III | ca. 3200-3000 |
| Early Dynastic I-II (ED I-II) | ca. 2900-2700 |
| Early Dynastic IIIa (ED IIIa = Fara period) | ca. 2600-2500 |
| Early Dynastic IIIb (ED IIIb) | ca. 2500-2340 |
| Old Akkadian (Sargonic) | ca. 2340-2200 |
| Third Dynasty of Ur (Ur III) | ca. 2100-2000 |
| Old Babylonian (OB) | ca. 2000-1600 |

# Appendix 9.C: Transliteration and translation of text 5 (CUNES 50-08-001)

Obverse, col. *i*

**(Table A)**

| | **Transliteration** | **Translation** |
|---|---|---|
| 1. | 1(aš) ninda-DU sa$_2$ | 1 *ninda* equal side |
| 2. | 1(aš) sar | 1 *sar* |
| 3. | 2(aš) sa$_2$ | 2 *ninda* equal side |
| 4. | 4(aš) sar | 4 *sar* |
| 5. | 3(aš) sa$_2$ | 3 *ninda* equal side |
| 6. | 1(u) la$_2$ 1(aš) sar | 10 minus 1 *sar* |
| 7. | 4(aš) sa$_2$ | 4 *ninda* equal side |
| 8. | 1(u) 6(aš) sar | 16 *sar* |
| 9. | 5(aš) sa$_2$ | 5 *ninda* equal side |
| 10. | 2(u) 5(aš) sar | 25 *sar* |
| 11. | 6(aš) sa$_2$ | 6 *ninda* equal side |
| 12. | 3(u) 6(aš) sar | 36 *sar* |
| 13. | 7(aš) sa$_2$ | 7 *ninda* equal side |
| 14. | 5(u) la$_2$ 1(aš) sar | 50 minus 1 *sar* |
| 15. | 8(aš) sa$_2$ | 8 *ninda* equal side |
| 16. | 1(geš$_2$) 4(aš) sar | 1 *geš* 4 *sar* |
| 17. | 9(aš) sa$_2$ | 9 *ninda* equal side |
| 18. | 1(geš$_2$) 2(u) 1(aš) sar | 1 *geš* 21 *sar* |
| 19. | 1(u) sa$_2$ | 10 *ninda* equal side |
| 20. | 1(iku) GAN$_2$ | 1 *iku GAN* |
| 21. | 2(u) sa$_2$ | 20 *ninda* equal side |
| 22. | 4(iku) GAN$_2$ | 4 *iku GAN* |

Obverse, col. *ii*

| | | |
|---|---|---|
| 1. | 3(u) sa$_2$ | 30 *ninda* equal side |
| 2. | 1(eše$_3$) 3(iku) GAN$_2$ | 1 *eše GAN* |
| 3. | 4(u) sa$_2$ | 40 *ninda* equal side |
| 4. | 2(eše$_3$) 4(iku) GAN$_2$ | 1 *eše* 4 *iku GAN* |
| 5. | 5(u) sa$_2$ | 50 *ninda* equal side |
| 6. | 1(bur$_3$) 1(eše$_3$) 1(iku) GAN$_2$ | 1 *bur* 1 *eše* 1 *iku GAN* |
| 7. | 1(geš$_2$) sa$_2$ | 60 *ninda* equal side |
| 8. | 2(bur$_3$) GAN$_2$ | 2 *bur GAN* |
| 9. | 2(geš$_2$) sa$_2$ | 2×60 *ninda* equal side |
| 10. | 8(bur$_3$) GAN$_2$ | 8 *bur GAN* |
| 11. | 3(geš$_2$) sa$_2$ | 3×60 *ninda* equal side |
| 12. | 1(bur’u) 8(bur$_3$) GAN$_2$ | 1 *buru* 8 *bur GAN* |
| 13. | 4(geš$_2$) sa$_2$ | 4×60 *ninda* equal side |
| 14. | 3(bur’u) 2(bur$_3$) GAN$_2$ | 3 *buru* 2 *bur GAN* |
| 15. | 5(geš$_2$) sa$_2$ | 5×60 *ninda* equal side |
| 16. | 5(bur’u) GAN$_2$ | 5 *buru GAN* |
| 17. | 6(geš$_2$) sa$_2$ | 6×60 *ninda* equal side |
| 18. | 1(šar$_2$) 1(bur’u) 2(bur$_3$) GAN$_2$ | 1 *šar* 1 *buru* 2 *bur GAN* |
| 19. | 7(geš$_2$) sa$_2$ | 7×60 *ninda* equal side |
| 20. | 1(šar$_2$) 3(bur’u) 8(bur$_3$) GAN$_2$ | 1 *šar* 3 *buru* 8 *bur GAN* |
| 21. | 8(geš$_2$) sa$_2$ | 8×60 *ninda* equal side |
| 22. | 2(šar$_2$) 8(bur$_3$) GAN$_2$ | 1 *šar* 8 *bur GAN* |

Obverse, col. *iii*

| | | |
|---|---|---|
| 1. | 9(geš$_2$) sa$_2$ | 9×60 *ninda* equal side |
| 2. | 2(šar$_2$) 4(bur’u) 2(bur$_3$) GAN$_2$ | 2 *šar* 4 *buru* 2 *bur GAN* |
| 3. | 1(geš’u) sa$_2$ | 600 *ninda* equal side |
| 4. | 3(šar$_2$) 2(bur’u) GAN$_2$ | 3 *šar* 2 *buru GAN* |
| 5. | 2(geš’u) sa$_2$ | 2×600 *ninda* equal side |
| 6. | 1(šar’u) 3(šar$_2$) 2(bur’u) GAN$_2$ | 1 *šaru* 1 *šar* 1 *buru GAN* |
| 7. | 3(geš'u) sa$_2$ | 3×600 *ninda* equal side |
| 8. | 3(šar'u) GAN$_2$ | 3 *šaru GAN* |
| 9. | 4(geš'u) sa$_2$ | 4×600 *ninda* equal side |
| 10. | 5(šar'u) 3(šar$_2$) 2(bur'u) GAN$_2$ | 5 *šaru* 3 *šar* 2 *buru GAN* |
| 11. | 5(geš'u) sa$_2$ | 5×600 *ninda* equal side |
| 12. | 1(šar$_2$) 2(šar'u) gal[36] 3(šar$_2$) 2(bur'u) GAN$_2$ | 1 large-*šar* 2 *šaru* 3 *šar* 2 *buru GAN* |
| 13. | 1(šar$_2$) sa$_2$ | 3600 *ninda* equal side |
| 14. | 2(šar$_2$) gal GAN$_2$ | 2 large-*šar GAN* |
| 15. | 2(šar$_2$) sa$_2$ | 2×3600 *ninda* equal side |
| 16. | 8(šar$_2$) gal GAN$_2$ | 8 large-*šar GAN* |
| 17. | 3(šar$_2$) sa$_2$ | 3×3600 *ninda* equal side |

Obverse, col. *iv*

| | | |
|---|---|---|
| 1. | 1(šar'u) 8(šar$_2$) gal GAN$_2$ | 2 large-*šaru* 1 large-*šar GAN* |
| 2. | 4(šar$_2$) sa$_2$ | 4×3600 *ninda* equal side |
| 3. | 3(šar'u) 2(šar$_2$) gal GAN$_2$ | 3 large-*šaru* 2 large-*šar GAN* |
| 4. | 5(šar$_2$) sa$_2$ | 5×3600 *ninda* equal side |
| 5. | 5(šar'u) gal GAN$_2$ | 5 large-*šaru GAN* |
| 6. | 6(šar$_2$) sa$_2$ | 6×3600 *ninda* equal side |
| 7. | 1(šar$_2$) KID[37] 1(šar'u) 2(šar$_2$) gal GAN$_2$ | 1 super-*šar* 1 large-*šaru* 2 large-*šar GAN* |
| 8. | 7(šar$_2$) sa$_2$ | 7×3600 *ninda* equal side |
| 9. | 1(šar$_2$) KID 3(šar'u) 8(šar$_2$) gal GAN$_2$ | 1 super-*šar* 3 large-*šaru* 8 large-*šar GAN* |
| 10. | 8(šar$_2$) sa$_2$ | 8×3600 *ninda* equal side |
| 11. | 2(šar$_2$) KID 8(šar$_2$) gal GAN$_2$ | 2 super-*šar* 8 large-*šar GAN* |
| 12. | 9(šar$_2$) sa$_2$ | 9×3600 *ninda* equal side |

Obverse, col.*v*

| | | |
|---|---|---|
| 1. | 2(šar$_2$) KID 4(šar'u) 2(šar$_2$) gal GAN$_2$ | 2 super-*šar* 4 large-*šaru* 2 large-*šar GAN* |
| 2. | 1(šar'u) sa$_2$ | 36000 *ninda* equal side |
| 3. | 3(šar$_2$) KID 2(šar'u) gal GAN$_2$ | 3 super-*šar* 2 large-*šaru GAN* |

**(Table B)**

| | | |
|---|---|---|
| 4. | 1(aš) nig$_2$-kas$_7$ sa$_2$ | 1 *nikkas* equal side |
| 5. | 3(aš) 2/3 5 gin$_2$ | 3 2/3 *gin* 5 (*gin-bi*) |
| 6. | 2(aš) nig$_2$-kas$_7$ sa$_2$ | 2 *nikkas* equal side |
| 7. | {sar} 15 gin$_2$ | 15 *gin* |
| 8. | 3(aš) nig$_2$-kas$_7$ sa$_2$ | 3 *nikkas* equal side |
| 9. | ½ sar 3(aš) 2/3 5 gin$_2$ | ½ *sar* 3 2/3 *gin* 5 (*gin-bi*) |
| 10. | 4(aš) nig$_2$-kas$_7$ sa$_2$ | 4 *nikkas* equal side |
| 11. | 1(aš) sar | 1 *sar* |
| 12. | 5(aš) nig$_2$-kas$_7$ sa$_2$ | 5 *nikkas* equal side |
| 13. | 1(aš) ½ sar 3(aš) 2/3 5 gin$_2$ | 1 ½ *sar* 3 2/3 *gin* 5 (*gin-bi*) |
| 14. | 6(aš) nig$_2$-kas$_7$ sa$_2$ | 6 *nikkas* equal side |
| 15. | sar 2(aš)[38] 15 gin$_2$ | 2 *sar* 15 *gin* |
| 16. | 7(aš) nig$_2$-kas$_7$ sa$_2$ | 7 *nikkas* equal side |
| 17. | 3(aš) sar 3(aš) 2/3 5 gin$_2$ | 3 *sar* 3 2/3 *gin* 5 (*gin-bi*) |

Obverse, col. *vi*

| | | |
|---|---|---|
| 1. | 8(aš) nig$_2$-kas$_7$ sa$_2$ | 8 *nikkas* equal side |
| 2. | 4(aš) sar | 4 *sar* |
| 3. | 9(aš) nig$_2$-kas$_7$ | 9 *nikkas* equal side |
| 4. | 5(aš) sar 3(aš) 2/3 5 gin$_2$ | 5 *sar* 3 2/3 *gin* 5 (*gin-bi*) |
| 5. | 1(u) nig$_2$-kas$_7$ sa$_2$ | 10 *nikkas* equal side |
| 6. | 6(aš) sar 15 gin$_2$ | 6 *sar* 15 *gin* |

---

36 The order of the signs is '1(šar$_2$) 2(šar'u) gal'; however, it is clear that the qualifier 'gal' refers to '1(šar$_2$)' and not to '2(šar'u)'. '*N* gal' means '60 times *N*'.

37 '*N* KID' means '60×60 times *N*' (Friberg 2007: 420).

38 The order is reversed: 'sar 2(aš)' is noted instead of '2(aš) sar'.

## (Table C)

| | | |
|---|---|---|
| 7. | 1(aš) kuš$_3$-numun sa$_2$ | 1 *kuš-numun* equal side |
| 8. | 1(aš) 2/3 gin$_2$ | 1 2/3 *gin* |
| 9. | 2(aš) kuš$_3$-numun sa$_2$ | 2 *kuš-numun* equal side |
| 10. | 6(aš) 2/3 gin$_2$ | 6 2/3 *gin* |
| 11. | 3(aš) kuš$_3$-numun sa$_2$ | 3 *kuš-numun* equal side |
| 12. | {sar} 15 gin$_2$ | 15 *gin* |
| 13. | 4(aš) kuš$_3$-numun sa$_2$ | 4 *kuš-numun* equal side |
| 14. | sar 1/3[39] 6(aš) 2/3 gin$_2$ | 1/3 *sar* 6 2/3 *gin* |
| 15. | 5(aš) kuš$_3$-numun sa$_2$ | 5 *kuš-numun* equal side |
| 16. | sar 2/3[40] 1(aš) 2/3 gin$_2$ | 2/3 *sar* 1 2/3 *gin* |
| 17. | 6(aš) kuš$_3$-numun sa$_2$ | 6 *kuš-numun* equal side |
| 18. | 1(aš) sar | 1 *sar* |
| 19. | 7(aš) kuš$_3$-numun sa$_2$ | 7 *kuš-numun* equal side |

Obverse, col. *vii*

| | | |
|---|---|---|
| 1. | 1(aš) sar 1/3 1(aš) 2/3 gin$_2$ | 1 1/3 *sar* 1 2/3 *gin* |
| 2. | 8(aš) kuš$_3$-numun sa$_2$ | 8 *kuš-numun* equal side |
| 3. | 1(aš) 2/3 sar 6(aš) 2/3 gin$_2$ | 1 2/3 *sar* 6 2/3 *gin* |
| 4. | 9(aš) kuš$_3$-numun sa$_2$ | 9 *kuš-numun* equal side |
| 5. | [2(aš)] sar 15 gin$_2$ | 2 *sar* 15 *gin* |
| 6. | 1(u) kuš$_3$-numun sa$_2$ | 10 *kuš-numun* equal side |
| 7. | 2(aš) 2/3 sar 6(aš) 2/3 gin$_2$ | 2 2/3 *sar* 6 2/3 *gin* |

## (Table D)

| | | |
|---|---|---|
| 8. | 1(aš) giš-bad sa$_2$ | 1 *giš-bad* equal side |
| 9. | gin$_2$ 1/3[41] 5(aš) gin$_2$-bi | 1/3 *gin* 5 *gin-bi* |
| 10. | 2(aš) giš-bad sa$_2$ | 2 *giš-bad* equal side |
| 11. | 1 2/3 gin$_2$ | 1 2/3 *gin* |
| 12. | 3(aš) giš-bad sa$_2$ | 3 *giš-bad* equal side |
| 13. | 3 2/3 gin$_2$ 5(aš) gin$_2$-bi | 3 2/3 *gin* 5 *gin-bi* |
| 14. | 4(aš) giš-bad sa$_2$ | 4 *giš-bad* equal side |
| 15. | 6 2/3 gin$_2$ | 6 2/3 *gin* |
| 16. | 5(aš) giš-bad sa$_2$ | 5 *giš-bad* equal side |
| 17. | {sar} 10 1/3 5 gin$_2$-bi gin$_2$[42] | 10 1/3 *gin* 5 *gin-bi* |
| 18. | 6(aš) giš-bad sa$_2$ | 6 *giš-bad* equal side |

Reverse, col. *i*

| | | |
|---|---|---|
| 1. | {sar}15 gin$_2$ | 15 *gin* |
| 2. | 7(aš) giš-bad sa$_2$ | 7 *giš-bad* equal side |
| 3. | sar 1/3 1/3 5 gin$_2$-bi gin$_2$[43] | 1/3 *sar* 1/3 *gin* 5 *gin-bi* |
| 4. | 8(aš) giš-bad sa$_2$ | 8 *giš-bad* equal side |
| 5. | sar 1/3 6 2/3 gin$_2$[44] | 1/3 *sar* 6 2/3 *gin* |
| 6. | 9(aš) giš-bad sa$_2$ | 9 *giš-bad* equal side |
| 7. | sar 1/2[45] 3 2/3 gin$_2$ 5 gin$_2$-bi | 1/2 *sar* 3 2/3 *gin* 5 *gin-bi* |
| 8. | 1(u) giš-bad sa$_2$ | 10 *giš-bad* equal side |
| 9. | sar 2/3[46] 1 2/3 gin$_2$ | 2/3 *sar* 1 2/3 *gin* |

---

39 The order is reversed: 'sar 1/3' is noted instead of '1/3 sar'.

40 The order is reversed: 'sar 2/3' is noted instead of '2/3 sar'.

41 The order is reversed: 'gin$_2$ 1/3' is noted instead of '1/3 gin$_2$'.

42 Order of the signs unclear.

43 Order of the signs unclear.

44 Order of the signs unclear.

45 The order is reversed: 'sar 1/2' is noted instead of '1/2 sar'.

46 The order is reversed: 'sar 1/2' is noted instead of '1/2 sar'.

**(Table E)**

| | | |
|---|---|---|
| 10. | 1(aš) šu-bad sa2 | 1 *šu-bad* equal side |
| 11. | gin2-bi-TA? 6(aš) 1(u) 5(aš) gin2-ba-gin2[47] | 6 *gin-bi* 15 *gin-ba-gin* |
| 12. | 2(aš) šu-bad sa2 | 2 *šu-bad* equal side |
| 13. | gin2 1/3[48] 5(aš) gin2-bi | 1/3 *gin* 5 *gin-bi* |
| 14. | 3(aš) šu-bad sa2 | 3 *šu-bad* equal side |
| 15. | 5(u) 6(aš) gin2-(bi) 15 gin2-ba-gin2 | 56 *gin-bi* 15 *gin-ba-gin*[49] |

Reverse, col. *ii*

| | | |
|---|---|---|
| 1. | 4(aš) šu-bad sa2 | 4 *šu-bad* equal side |
| 2. | 1 2/3 gin2 | 1 2/3 *gin* |
| 3. | 5(aš) šu-bad sa2 | 5 *šu-bad* equal side |
| 4. | 2 ½ gin2 5[sic](aš) gin2-TA?-bi 15 ba-gin2-gin2[50] | 2 1/2 *gin* 6! *gin-bi* 15 *gin-ba-gin* |
| 5. | 6(aš) šu-bad sa2 | 6 *šu-bad* equal side |
| 6. | 3 2/3 gin2 5(aš) gin2-bi | 3 2/3 *gin* 5 *gin-bi* |
| 7. | 7(aš) šu-bad sa2 | 7 *šu-bad* equal side |
| 8. | 5(aš) gin2 gin2-bi 6(aš) 15 gin2-ba-gin2 | 5 *gin* 6 *gin-bi* 15 *gin-ba-gin* |
| 9. | 8(aš) šu-bad sa2 | 8 *šu-bad* equal side |
| 10. | 6 2/3 gin2 | 6 2/3 *gin* |
| 11. | 9(aš) šu-bad sa2 | 9 *šu-bad* equal side |
| 12. | 8 1/3 gin2 gin2-bi-TA 6(aš) 1(u) 5(aš) gin2-ba-gin2-TA | 8 1/3 *gin* 6 *gin-bi* 15 *gin-ba-gin* |

Reverse, col. *iii*

| | | |
|---|---|---|
| 1. | 1(u) šu-bad sa2 | 10 *šu-bad* equal side |
| 2. | {sar} 10 1/3 gin2 5 gin2-bi | 10 1/3 gin 5 gin-bi[51] |

**(subscript - transliteration / translation CDLI, 2012-10-12 by Englund)**

| | | |
|---|---|---|
| 3. | ka9-ka9# | ... |
| 4. | ka9 ku3# | ... |
| 5. | lu2 dub-sar | *scribe* |
| 6. | lugal-he2-gal2-su3 | *Lugal-ḫegal-su* |
| 7. | er-da | *Erda,* |
| 8. | sanga# ku3-ma2? | ... |
| 9. | sanga ka9-še3 mu-zu-kur? | … |

47 The order of the signs is unclear.

48 The order is reversed: ‘gin2 1/3’ is noted instead of ‘1/3 gin2’.

49 The scribes noted ‘5(u) 6 (aš) gin2-(bi)’ (56 *gin-bi*) instead of the expected equivalent notation using fractions of *gin* ‘2/3 gin2 16 gin2-bi’ (2/3 *gin* 16 *gin-bi*).

50 Order of the signs unclear.

51 The order of the signs is slightly different from other instance of the same surface (in Table D), which may indicate that the calculation was executed again, and not copied from the previous table.

# References

## Abbreviations

| | |
|---|---|
| A | Tablets in the collections of the Oriental Institute, Univ. of Chicago |
| BCE | Before the Common Era |
| Col. | Column |
| CUNES | Tablets in the collection of the Cornell University Near Eastern Studies |
| MS | Tablets in the Schøyen collection |
| Obv. | Obverse |
| Rev. | Reverse |
| Translit. | Transliteration |
| VAT | Tablets in the *Vorderasiatische Abteilung Tontafeln*, Berlin |

## List of primary Sources

| | | |
|---|---|---|
| VAT 12593 (text 1) | Deimel 1923: No 82 (SF 82) | P010678 |
| MS 3047 (text 2) | Friberg 2007: 160 | P252059 |
| Feliu 2012 (text 3) | Feliu 2012 | P464229 |
| A 681 (text 4) | Luckenbill 1930: 36 | P222256 |
| CUNES 50-08-001 (text 5) | Friberg 2007: 419-426. | P274845 |

## Bibliography

Biggs, Robert D. 1974. *Inscriptions from Tell Abu Salabikh, Oriental Institute Publications.* Chicago: The University of Chicago Press.

Chemla, Karine. Forthcoming. Tables as texts. In *The history of numerical tables*, ed. Dominique Tournès, Chap. 1. Berlin: Springer.

Colonna d'Istria, Laurent. 2015. La notation des fractions dans la documentation de l'époque des derniers *šakkanakkū* de Mari. *Akkadica* 136: 103-125.

Damerow, Peter. 2016. The Impact of Notation Systems: From the Practical Knowledge of Surveyors to Babylonian Geometry. In *Spatial Thinking and External Representation. Towards a Historical Epistemology of Space*, ed. Matthias Schemmel, 93-119. Berlin: Max Planck Research Library for the History and Development of Knowledge. Edition Open Access.

Deimel, Anton. 1923. *Die Inschriften von Fara II. Schultexte aus Fara*. Vol. 43, *Wissenschaftliche Veröffentlichungen der Deutschen Oriental Gesellschaft*. Leipzig: J. C. Hinrichs'sche Buchhandlung

Edzard, Dietz Otto. 1969. Eine altsumerische Rechentafel. In *lišan mithurti, Festschrift Wolfram Freiherr von Soden*, ed. Wolfgang Röllig, 101-104. Neukirchen-Vluyn: Neukirchener Verlag des Erziehungs-vereins.

Feliu, Lluís. 2012. A new Early Dynastic IIIb metro-mathematical table tablet of area measures from Zabalam. *Altorientalische Forschungen* 39:218-225.

Friberg, Jöran. 1987-1990. Mathematik. *Reallexikon der Assyriologie* 7:531-585.

Friberg, Jöran. 1997-1998. Round and almost round numbers in proto-literate metro-mathematical field texts. *Archiv für Orientforshung* 44-45: 3-58.

Friberg, Jöran. 2007. *A Remarkable Collection of Babylonian Mathematical Texts*. Vol. I, *Manuscripts in the Schøyen Collection: Cuneiform Texts*. New York: Springer.

Friberg, Jöran. 2019. Three thousand years of sexagesimal numbers in Mesopotamian mathematical texts. *Archive for History of Exact Sciences* 73: 183–216.

Høyrup, Jens. 2002. *Lengths, Widths, Surfaces. A Portrait of Old Babylonian Algebra and its Kin, Studies and Sources in the History of Mathematics and Physical Sciences*. Berlin, London: Springer.

Jestin, Raymond R. 1937. *Tablettes sumériennes de Šuruppak conservées au musée de Stamboul.* Paris: de Boccard.

Jestin, Raymond R . 1957. *Nouvelles tablettes de Šuruppak du musée de Stamboul*. Paris: de Boccard.

Krebernik, Manfred. 1998. Die Texte aus Fara und Tell Abu Salabih. In *Mesopotamien. Späturuk-Zeit und Frühdynastische Zeit*, ed. Josef Bauer, Robert K. Englund and Manfred Krebernik, 235-427. Freiburg, Göttingen: Universitätsverlag.

Liverani, Mario. 1996. Reconstructing the rural landscape of the Ancient Near East. *Journal of the Economic and Social History of the Orient* 39 (1):1-41.

Luckenbill, Daniel D. 1930. *Inscriptions from Adab,*. Chicago: Oriental Institute Publications.

Melville, Duncan J. 2002. Ration computations at Fara: multiplication or repeated addition. In *Under One Sky. Astronomy and Mathematics in the Ancient Near East*, ed. John. M. Steele and Annette Imhausen, 237-252. Münster: Ugarit-Verlag.

Neugebauer, Otto. 1935. *Mathematische Keilschrifttexte I.* Berlin: Springer.

Nissen, Hans J. , Peter Damerow and Robert Englund. 1993. *Archaic Bookkeeping. Writing and Techniques of Economic Administration in the Ancient Near East.* Chicago: University of Chicago Press.

Ossendrijver, Mathieu, Christine Proust, Micah Ross and Nathan Sidoli. Forthcoming. Ancient Near-East tables (Egypt, Mesopotamia, Greece). In *The History of Numerical Tables*, ed. Dominique Tournès. Berlin: Springer.

Ouyang Xiaoli and Christine Proust. forthcoming. Place value notations in the Ur III period: marginal numbers in administrative texts. In *Cultures of Computation and Quantification*, ed. Karine Chemla, Agathe Keller and Christine Proust, Heidelberg, New York: Springer.

Powell, Marvin A. 1972. Sumerian area measures and the alleged decimal substratum. *Zeitschrift für Assyriologie und Vorderasiatische Archäologie* 62:165–222.

Powell, Marvin A. 1973. On the reading and meaning of GANA2. *Journal of Cuneiform Studies* 25: 178-184.

Powell, Marvin A. 1976. The antecedents of Old Babylonian place notation and the early history of Babylonian mathematics. *Historia Mathematica* 3:417-439.

Powell, Marvin A. 1987-1990. Masse und Gewichte. In *Reallexikon der Assyriologie*. Ed. Edzard, Dietz Otto. Berlin: De Gruyter.

Proust, Christine. 2009. Numerical and metrological graphemes: from cuneiform to transliteration. *Cuneiform Digital Library Journal.* http://www.cdli.ucla.edu/pubs/cdlj/2009/cdlj2009_001.html.

Robson, Eleanor. 2003. Tables and tabular formatting in Sumer, Babylonia, and Assyria, 2500 BCE - 50 CE. In *History of Mathematical Tables: from Sumer to Spreadsheets*, ed. Martin Campbell-Kelly, Mary Croarken, Raymond Flood and E. Robson. Oxford: Oxford University Press.

Tournès, Dominique (ed.). forthcoming. *The history of Numerical Tables*. Berlin: Springer.